\documentclass[3p,review]{elsarticle}

\usepackage{amssymb}
\usepackage{amsmath}
\usepackage{amsthm}

\usepackage{float}
\usepackage[pdfborder={0 0 0}, colorlinks=true, linkcolor=blue, citecolor=red]{hyperref}
\hypersetup{pdfborder=0 0 0}
\usepackage{algorithm}
\usepackage{algpseudocode}
\usepackage{mathrsfs}
\usepackage{subfigure}
\usepackage{multirow}
\usepackage{wrapfig}
\usepackage{siunitx} 

\usepackage{leftidx}
\usepackage[latin1]{inputenc}
\usepackage{ulem}

\usepackage{IEEEtrantools}

\usepackage{lineno,hyperref}
\modulolinenumbers[5]

\usepackage[figuresright]{rotating}

\def\e{{\epsilon}}

\def\norm#1{\|#1\|} 

\newcommand{\beq} {\begin{equation}}
\newcommand{\eeq} {\end{equation}}
\newcommand{\bdm} {\begin{displaymath}}
\newcommand{\edm} {\end{displaymath}}
\newcommand{\bit}{\begin{itemize}}
\newcommand{\eit}{\end{itemize}}
\newcommand{\bde}{\begin{description}}
\newcommand{\ede}{\end{description}}
\newcommand{\bce}{\begin{center}}
\newcommand{\ece}{\end{center}}
\newcommand{\ben} {\begin{enumerate}}
\newcommand{\een} {\end{enumerate}}
\newcommand{\bea} {\begin{eqnarray}}
\newcommand{\eea} {\end{eqnarray}}
\newcommand{\barr} {\begin{array}}
\newcommand{\earr} {\end{array}}
\newcommand{\bean} {\begin{eqnarray*}}
\newcommand{\eean} {\end{eqnarray*}}
\newcommand{\edoc} {\end{document}}

\newcommand*\ddiff{\mathop{}\!\mathrm{d}}

\newcommand{\dd}[2] {\ensuremath{\frac{\partial {#1}}{\partial {#2}}}}
\newcommand{\DD}[2] {\ensuremath{\frac{\ddiff {#1}}{\ddiff {#2}}}}

\newcommand{\eval}[2][\right]{\relax
  \ifx#1\right\relax \left.\fi#2#1\rvert}

\providecommand{\abs}[1]{\left\lvert#1\right\rvert}
\providecommand{\norm}[1]{\left\hspace{-0.5ex}\lVert#1\right\hspace{-0.5ex}\rVert}

\newcommand{\Nel} {\ensuremath{{N_\text{el}}}}

\newcommand{\half} {\ensuremath{\frac{1}{2}}}

\newcommand{\mc}[1]{\mathcal{#1}}

\newcommand{\mb}[1]{\mathbf{#1}}

\newcommand{\nor}[1]{\left\vert\kern-0.35ex\left\vert #1 \right\vert\kern-0.35ex\right\vert}
\newcommand{\norL}[1]{\left\vert\kern-0.5ex\left\vert #1 \right\vert\kern-0.5ex\right\vert}
\newcommand{\nort}[1]{{\left\vert\kern-0.35ex\left\vert\kern-0.35ex\left\vert #1 
    \right\vert\kern-0.35ex\right\vert\kern-0.35ex\right\vert}}

\newcommand{\LRp}[1]{\left( #1 \right)}
\newcommand{\LRs}[1]{\left[ #1 \right]}
\newcommand{\LRa}[1]{\left< #1 \right>}
\newcommand{\LRc}[1]{\left\{ #1 \right\}}

\newcommand{\jumpL}[1] {\ensuremath{\LRs{\hspace{-1ex}\LRs{#1}\hspace{-1ex}}}}

\newcommand{\average}[1] {\ensuremath{\LRc{\hspace{-0.9ex}\LRc{#1}\hspace{-0.9ex}}}}
\newcommand{\averageM}[1] {\ensuremath{\LRc{\hspace{-1.ex}\LRc{#1}\hspace{-1.ex}}}}

\newcounter{tablecell}

\newcommand{\figlab}[1]{\label{fig:#1}}
\newcommand{\eqnlab}[1]{\label{eq:#1}}

\newcommand{\tablab}[1]{\label{tab:#1}}

\newcommand{\figref}[1]{\ref{fig:#1}}

\newcommand{\eqnref}[1]{\eqref{eq:#1}}

\newcommand{\seclab}[1]{\label{sect:#1}}
\newcommand{\secref}[1]{\ref{sect:#1}}
\newcommand{\tabref}[1]{\ref{tab:#1}}

\renewcommand{\u}{u}

\newcommand{\q}{q}

\newcommand{\R}{{\mathbb{R}}}

\newcommand{\Poly}{{\mc{P}}}

\newcommand{\V}{{V}}

\newcommand{\Vhp}[1]{{\V_h^#1}}

\newcommand{\n}{\mb{n}}
\newcommand{\nm}{\n^-}
\newcommand{\np}{\n^+}

\renewcommand{\sb}{\mb{s}}

\newcommand{\Tmesh}{\mc{T}}
\newcommand{\pTmesh}{\partial\mc{T}}
\newcommand{\Tmeshh}{{{\Tmesh}_h}}

\newcommand{\nb}{{\bf n}}

\newcommand{\Fcal}{\mathcal{F}}

\newcommand{\dt}{{\triangle t}}
\newcommand{\dx}{{\triangle x}}

\newcommand{\Ical}{\mathcal{I}}

\newcommand{\Xhat}{{\hat{X}}}

\newcommand{\Yhat}{{\hat{Y}}}
\newcommand{\Zhat}{{\hat{Z}}}

\renewcommand{\u}{{u}}

\newcommand{\Elem}{{I}}

\usepackage{scalerel,stackengine}
\stackMath
\newcommand\reallywidehat[1]{%
\savestack{\tmpbox}{\stretchto{%
  \scaleto{%
    \scalerel*[\widthof{\ensuremath{#1}}]{\kern-.6pt\bigwedge\kern-.6pt}%
    {\rule[-\textheight/2]{1ex}{\textheight}}
  }{\textheight}%
}{0.5ex}}%
\stackon[1pt]{#1}{\tmpbox}%
}

\usepackage{color}
\usepackage{soul,xargs}
\usepackage[pdftex,dvipsnames]{xcolor}

\usepackage[colorinlistoftodos,prependcaption,textsize=tiny]{todonotes}

\newcommandx{\question}[2][1=]{\todo[linecolor=red,backgroundcolor=red!25,bordercolor=red,#1]{#2}}
\newcommandx{\change}[2][1=]{\todo[linecolor=blue,backgroundcolor=blue!25,bordercolor=blue,#1]{#2}}
\newcommandx{\add}[2][1=]{\todo[linecolor=OliveGreen,backgroundcolor=OliveGreen!25,bordercolor=OliveGreen,#1]{#2}}
\newcommandx{\improve}[2][1=]{\todo[linecolor=Plum,backgroundcolor=Plum!25,bordercolor=Plum,#1]{#2}}
\newcommandx{\thiswillnotshow}[2][1=]{\todo[disable,#1]{#2}}
\newcommandx{\remove}[2][1=]{\todo[linecolor=yelllow,backgroundcolor=yellow!10,bordercolor=red,#1]{#2}}

\begin{document}

\begin{frontmatter}




\title{Learning Subgrid-Scale Models with Neural Ordinary Differential Equations}

\author[AddrANL]{Shinhoo Kang\corref{mycorrespondingauthor}}
\cortext[mycorrespondingauthor]{Corresponding author}
\ead{shinhoo.kang@anl.gov}
\address[AddrANL]{Mathematics and Computer Science Division, Argonne National Laboratory. Lemont, IL 60439, USA.}

\author[AddrANL]{Emil M. Constantinescu} 
\ead{emconsta@anl.gov}
 

\begin{abstract} 
  
We propose a new approach to learning the subgrid-scale model when simulating partial differential equations (PDEs) solved by the method of lines and their representation in chaotic ordinary differential equations, based on neural ordinary differential equations (NODEs). Solving systems with fine temporal and spatial grid scales is an ongoing computational challenge, and closure models are generally difficult to tune.
Machine learning approaches have increased the accuracy and efficiency of computational fluid dynamics solvers. In this approach neural networks are used to learn the coarse- to fine-grid map, which can be viewed as subgrid-scale parameterization. We propose a strategy that uses the NODE and partial knowledge to learn the source dynamics at a continuous level.
Our method inherits the advantages of NODEs and can be used to parameterize subgrid scales, approximate coupling operators, and improve the efficiency of low-order solvers. Numerical results with the two-scale Lorenz 96 ODE, the convection-diffusion PDE, and the viscous Burgers' PDE are used to illustrate this approach. 
  
\end{abstract}


\end{frontmatter}




\section{Introduction}
Many problems in science and engineering are expressed in the form of ordinary differential equations (ODEs),
\begin{align}
    \DD{u(t)}{t} = R(t,u(t)), 
    \eqnlab{ode}
\end{align}
where $u(t)$ is the state variable and $R$ is the functional that captures the space discretization and describes the overall dynamics. Problems in fluid dynamics and  numerical weather prediction, which are based on Navier--Stokes equations, are captured by this formalism. 
In practice, however, resolving all scales in physical systems is intractable. Therefore, subgrid-scale parameterization is used to account for the gap between the real-world system and its digital representation and underresolved-scale dynamics, which in the case of numerical weather models can include rain, clouds, effects of vegetation, and radiation. Such terms are formulated in terms of the forcing/source $S$ in \eqnref{ode},  
\begin{align}
    \DD{u(t)}{t} = R(t,u(t)) + S(u).
    \eqnlab{ode-source}
\end{align}
Such parameterizations interact with resolved-scale dynamics and prove to be crucial for enhancing predictability \cite{bauer2015quiet}. However, parameterizing subgrid-scale physics is not a trivial task \cite{hong2012next,irrgang2021towards}. 
For example, understanding the dynamics of the vegetation that covers the Earth's surface, absorbs solar radiation, and interacts with the atmosphere requires significant  research effort to explain its macroscopic effects on the Earth's scale \cite{franklin2020organizing}. Cloud parameterization is one of the major sources of  uncertainty in climate models \cite{schneider2017climate}.

Neural network parameterizations are particularly suitable to address the subgrid parametrization issues by virtue of their universal approximation  \cite{hornik1990universal}. 
Numerous studies have focused  on using finite-depth neural networks to match the source term $S(u)$ in \eqnref{ode-source} to data. 
Examples relevant to our work include  strategies such as that of Bolton and Zanna, who employed a convolution neural network to forecast unresolved small-scale turbulence~\cite{bolton2019applications}; 
Brenowitz and Bretheron used a feed-forward neural network to estimate the residual heating and moistening \cite{brenowitz2019spatially}; 
Rasp et al. investigated superparameterization of a cloud-resolving model using deep learning neural networks~\cite{rasp2018deep} and developed a machine learning parameterization~\cite{rasp2020coupled} for the coupling term in the two-scale Lorenz 96 model \cite{lorenz1996predictability}.
 Lara and Ferrer proposed  enhancing a low-order solver by introducing a parameterized source term~\cite{de2022accelerating,manrique2022accelerating}. 
The source term, according to the authors, acts as a corrective force to fill in the missing scales and interactions between low-order and high-order solutions. 
To that end, they were able to speed up the simulation and improve its accuracy by employing a feed-forward neural network to approximate the corrective forcing term. 
Huang et al. introduced a neural-network-based corrector for time integration methods to compensate for the temporal discretization error caused by coarse timestep size~\cite{huang2022accelerating}.
All these approaches involve learning a static correction, however, which may make the simulation sensitive to the discretization parameters. We propose a strategy that alleviates such a restriction by learning a continuous operator and representing its dynamics via neural ordinary differential equations (NODEs).

NODEs \cite{chen2018neural} are an elegant approach to learn the dynamics from time series data. This approach can be seen as residual neural networks (ResNets) \cite{he2016deep} in the continuous limit \cite{avelin2021neural}. That is, a NODE has an infinite number of layers that share the network parameters. This implies 
that the number of parameters in a NODE is smaller than that of ResNet \cite{gupta2022galaxy}; hence, a NODE can also be memory efficient.  
A NODE can choose any subinterval in time so long as the integration is stable, and therefore it can handle nonuniform samples or missing data in time series 
\cite{rubanova2019latent}. Also, various standard ODE solvers such as Runge--Kutta or Dormand--Prince methods \cite{dormand1980family} can be incorporated into a NODE. 
More important, the computationally efficient {\it adjoint sensitivity method} \cite{pontryagin1987mathematical} can be used to calculate gradients during learning.
This framework is powerful, and many promising results are reported 
in reversible generative flow~\cite{grathwohl2018ffjord}, 
image classification~\cite{zhuang2020adaptive,zhuang2021mali}, 
continuous normalizing flow~\cite{nguyen2019infocnf},
time series analysis~\cite{kidger2020neural}, and manifold learning~\cite{lou2020neural}. Recently, Rackauckas et al. proposed universal differential equations that embed a neural network in governing equations to utilize known physical information \cite{rackauckas2020universal}. 
Shankar et al. investigated machine-learning-based closure models for Burgers turbulence within a differentiable physics paradigm \cite{shankar2023differentiable}.

In this study we leverage NODEs and simulation data with different scales to learn the dynamics of the source term at a continuous level. 
Our approach inherits the benefits of NODEs and has the potential to be applicable to various parameterization problems. In this study, however, we do not take into account practical problems such as turbulent flows. Instead, we focus on simple problems for demonstrating our idea clearly. We will continue to explore our proposed methodology in more complex problems in the future. 
The contributions of our work are summarized as follows:
(1) we propose a novel approach to learn continuous source dynamics through neural ordinary differential equations; (2) we numerically demonstrate that the neural network source operator is capable of approximating coupling terms as well as reconstructing the subgrid-scale dynamics; and  
(3) we generalize the discrete corrective forcing approach \cite{de2022accelerating} to a continuous corrective forcing approach that supports high-order accuracy in time, and we compare the two methods.

The paper is organized as follows. 
We present the proposed method and model problems in Section \secref{nodeparam}.
In Section \secref{NumericalResults} several numerical examples are used to test our proposed methodology. 
We parameterize the coupling term between the slow and fast variables in the two-scale Lorenz 96 equation and the corrective forcing component in the convection-diffusion equation by using the proposed NODE approach. An example of Burgers turbulence with continuous corrective forcing is also investigated. 
In Section \secref{Conclusion} we present our conclusions.

\section{Neural ODE for Parameterizing the Source Function}
\seclab{nodeparam}

This section briefly reviews NODEs and provides the methodology for parameterizing the source term in \eqnref{ode-source}. 
Our presentation focuses on three problems: the two-scale Lorenz 96 system, the linear convection-diffusion equation, and the viscous Burgers' equation, which are utilized to approximate the coupling term and the corrective forcing terms, respectively. We also fix the neural network architecture to a configuration that has sufficient expressiveness, and we focus on the methodological aspects of the proposed strategy.

NODE \citep{chen2018neural} is a neural network model that uses ordinary differential equations to learn the dynamics of a stream of input-output data. 
In NODE, the functional $R$ in \eqnref{ode} is approximated by a neural network $R_\theta$ at a continuous level, 
\begin{align}
\eqnlab{node0}
    \DD{u(t)}{t} = R_\theta(t,u(t)).
\end{align}
Here, $R_\theta$ is parameterized by network weights and biases $\theta$, and $u(t)$ represents the state of NODE.
By integrating \eqnref{node0}, we can evaluate the state $u$ at any target time $T$, 
\begin{align}
\eqnlab{node0-integral}
    u(T) &= u(0) + \int_{0}^{T} R_\theta(\tau,u(\tau)) d\tau. 
\end{align}
We let $N_t$ be the number of subintervals from $0$ to $T$.
With $T_0=0$ and $T_{N_t}=T$, \eqnref{node0-integral} can be written as
\begin{align}
\eqnlab{node0-integral-Nt}
    u(T) &= u(0) + \sum_{n=1}^{N_t}\int_{T_{n-1}}^{T_n} R_\theta(\tau,u(\tau)) d\tau. 
\end{align}
We note that \eqnref{node0-integral-Nt} has only one neural network architecture. 
NODE does not require storing intermediate network parameters between subintervals. 
Thus, NODE can be more memory efficient than a traditional neural network such as ResNet \cite{he2016deep}, where $N_t$ network architectures are needed for $u(T)$.
NODE can also handle nonuniform samples by flexibly selecting any subinterval in time. 
Furthermore, backpropagation \cite{paszke2019pytorch} or adjoint methods \cite{chen2018neural} through the ODE solver can be used to train NODEs, allowing the neural network to modify its parameters and architecture according to the input data. This makes NODEs particularly useful for modeling time series data.

While NODEs have several benefits over traditional networks, they also have limitations. 
NODEs can be computationally expensive to train for complex models with many network parameters \cite{grathwohl2018ffjord,finlay2020train}. 
NODEs lack robustness and uncertainty modeling capabilities \cite{anumasa2021improving}. 
NODEs have no direct mechanism for incorporating data that arrives later 
\cite{kidger2020neural}. 
To address these issues, researchers have presented several variants of NODEs \cite{djeumou2022taylor,anumasa2021improving,kidger2020neural}. 
Nonetheless, we concentrate on the standard NODEs in our study as proof of concept.

Let us consider a neural-network-based source term $S_\theta(u)$ that approximates the source functional $S(u)$ in \eqnref{ode-source}: 
\begin{align}
    \DD{u(t)}{t} = R(t,u(t)) + S_\theta(u(t)).
    \eqnlab{ode-nnsource}
\end{align}
The neural network architecture is composed of an input layer, three internal layers, and an output layer, as shown in Figure 
\figref{neural-network-architecture}.
The input state is linearly connected to the input layer. Following that,  
the internal layers with $128$ neurons are connected by using rectified linear unit (ReLU) \cite{nair2010rectified} activation functions. 
The linear output layer yields
the source term $S_\theta$ evaluated on state $u$. 

The corrected solution is obtained by integrating \eqnref{ode-nnsource}.  We can thus evaluate the state $u$ at any target time ($t=T$) by 
\begin{align}
\eqnlab{node-integral}
    u(T) &= u(0) + \int_{0}^{T}  R(\tau,u(\tau)) + 
    S_\theta(u(\tau)) \ddiff\tau. 
\end{align}
Various timestepping methods can be utilized to approximate the integral of \eqnref{node-integral}. In this study, however, we  focus only on $s$-stage explicit Runge--Kutta (ERK) schemes with uniform timestep size.
Applying ERK schemes to \eqnref{node-integral} yields
\begin{align*}
  U_{n,i} &= u_n + \dt\sum_{j=1}^{i-1} a_{ij} \LRp{R_j + S_{\theta_j} }, \quad i=1,2,\cdots,s,\\
  u_{n+1} &= u_n + \dt\sum_{i=1}^{s}b_i \LRp{R_i+S_{\theta_i}},
\end{align*}
where $R_i:=R(t_n + c_i\dt,U_{n,i})$;
$S_{\theta_i}:=S_\theta(U_{n,i})$; and $a_{ij}$, $b_i$, and $c_i$ are scalar coefficients for $s$-stage ERK methods. Given an initial condition $u_0$, the $N$ step ERK timestepping method generates $N$ discrete instances of a single trajectory, $\LRc{u_1,u_2,\cdots,u_N}$. We denote the approximated solution $u$  of \eqnref{ode-nnsource} by $\hat{u}$ from now on.

\begin{wrapfigure}{r}{0.4\textwidth}
    \begin{center}
  \includegraphics[trim=3cm 0.8cm 4cm 5cm,clip=true,width=0.4\textwidth]{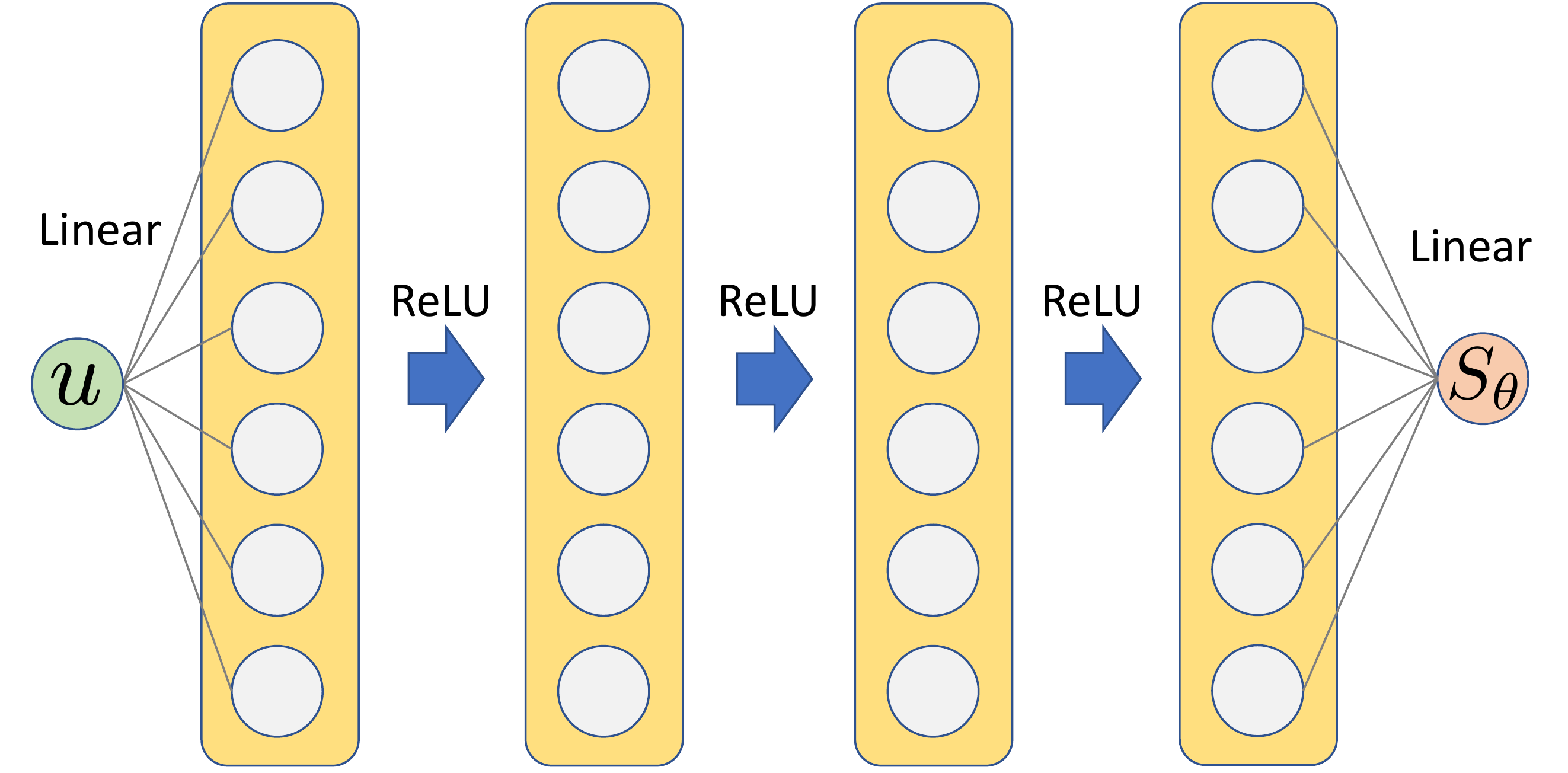}
  \end{center}
  \caption{Schematic of a neural network source.
  The input state $u$ is connected to the first layer. 
  The internal layers are then coupled together by using ReLU activation functions.
  The source term $S_\theta$ is produced from the linear output layer.
  }
  \figlab{neural-network-architecture}
\end{wrapfigure}

\subsection{Two-Scale Lorenz 96 Model}
The two-scale Lorenz 96 model is a simplified example of atmospheric dynamics with the terms of convection, diffusion, coupling, and forcing~ \citep{lorenz1996predictability}. The model dynamics is governed by two sets of  variables: the slow variables $X_k$, with $k=1,2,\cdots,K$, and the fast variables $Y_{j,k} (j=1,2,\cdots,J)$. 
The slow variables refer to atmospheric waves, whereas the fast variables stand for a faster dynamics over smaller scales such as convection.
The Lorenz 96 model has $N_{dof}=K(1+J)$ degrees of freedom, where slow and fast variables are coupled via the following system: 
\begin{subequations}
  \eqnlab{eq-lorenz96}
\begin{align}
  \eqnlab{eq-lorenz96-slow}
  \DD{X_k}{t} & = - X_{k-1} \LRp{X_{k-2}-X_{k+1} } - X_k + F - h \overline{Y}_k,\\
  \eqnlab{eq-lorenz96-fast}
  \frac{1}{c}\DD{Y_{j,k}}{t} &= - J Y_{j+1,k} \LRp{Y_{j+2,k} - Y_{j-1,k} } - Y_{j,k} + \frac{h}{J} X_k,
\end{align}
\end{subequations}
where $\overline{Y}_k :=\frac{1}{J}\sum_{j=1}^J Y_{j,k}$ is the average contribution of the fast variables $Y_{j,k}$ ($j=1,2,\cdots,J$) for a given $k$ corresponding to each slow variable; $c$ is the time-scale separation between $X_k$ and $Y_{j,k}$; $J$ controls the relative amplitude of the $Y_{j,k}$ component; and $h$ is the strength of the coupling between the fast and the slow variables \cite{carlu2019lyapunov}. Both $X_k$ and $Y_{j,k}$ have periodic boundary conditions, which represent a circle of constant latitude on Earth for the slow variables, whereas the fast ones represent internal grid variability or subgrid effects. 
Following \cite{lorenz1996predictability}, we choose $h=1$, $K=36$, and $c=J=10$. In this case, the fast variable is 10 times faster and smaller than the slow variable, and the dynamics is chaotic. 

The coupling term $-h \overline{Y}_k$ in \eqnref{eq-lorenz96-slow}  
is the contribution from the localized fast variables  
that are coupled with one slow variable $X_k$, and this proves to be an ideal setup to test the proposed methodology. To that end, 
we replace the coupling term with a neural approximation $\mc{S}_\theta (\Xhat_k)$, which yields 
  \begin{subequations}
  \eqnlab{eq-lorenz96-nnsource}
\begin{align}
  \eqnlab{eq-lorenz96-slow-nnsource}
  \DD{\Xhat_k}{t} & = - \Xhat_{k-1} \LRp{\Xhat_{k-2}-\Xhat_{k+1} } - \Xhat_k + F +\mc{S}_\theta(\Xhat_k),\\ 
  \eqnlab{eq-lorenz96-fast-nnsource}
  \frac{1}{c}\DD{\Yhat_{j,k}}{t} &= - J \Yhat_{j+1,k} \LRp{\Yhat_{j+2,k} - \Yhat_{j-1,k} } - \Yhat_{j,k} + \frac{h}{J} \Xhat_k.
\end{align}
\end{subequations}
This allows us to decouple the system \eqnref{eq-lorenz96}. Hence, 
we can solve the slow equation \eqnref{eq-lorenz96-slow-nnsource} with a larger stepsize compared with that of the coupled system \eqnref{eq-lorenz96-slow}. 

Using the fourth-order ERK (ERK4), we create several trajectories with random initial conditions for training data. A trajectory of the slow and fast variables is defined as $Z:=(X,Y)$. For training, 
we select $n$-batch trajectories and sampling time level $t_s$ at random and gather $m$ successive instances from each trajectory, that is, 
$\LRc{Z_{t_s}^{(i)},Z_{t_s+\dt}^{(i)},\cdots,Z_{t_s+m\dt}^{(i)} }$ for $i=1,\cdots,n$. 
Then, 
we predict the slow and fast variables 
$\LRc{\Zhat_{t_s+\dt}^{(i)},\Zhat_{t_s+2\dt}^{(i)},\cdots,\Zhat_{t_s+m\dt}^{(i)} }$ 
by solving \eqnref{eq-lorenz96-nnsource} given an initial condition $Z_{t_s}^{(i)}$ for $i=1,\cdots,n$. We use the mean squared error as the loss function:
\begin{align}
\eqnlab{lorenz96-loss}
  \mathrm{loss} =  \frac{1}{n m}\sum_{i=1}^n \sum_{\ell=1}^m \norm{Z_{t_s+\ell \dt}^{(i)} - \hat{Z}_{t_s+\ell \dt}^{(i)}  }^2. 
\end{align}
Here, $\norm{\cdot}$ is the $\ell^2$-norm. For example, $\norm{x}=\LRp{\sum_{i=1}^l x_i^2}^\half$ for a vector $x\in \R^l$. A common choice in many machine learning studies is the $\ell_2$ norm. Other norms can also be considered for computing loss functions. However, when we changed the loss function to the mean absolute error ($\ell_1$-norm), the prediction of the augmented DG approximation is not significantly different from that with the $\ell_2$-norm. Indeed, the authors in \cite{linot2023stabilized} also reported that the choice of $\ell_1$ and $\ell_2$ did not lead to substantial difference in their result for Burgers' equation. In this reason, we consistently used $\ell_2$ norm for the loss function in this study.

\subsection{Convection-Diffusion Equation} 
Our second example is the linear convection-diffusion equation defined on the time and space interval $(t,x) \in [0,T]\times\Omega$:
\begin{align}
  \eqnlab{pde-convdiff-eq}
 \dd{u}{t} + a\dd{u}{x} 
 = \kappa \dd{u}{x^2} \text{ in } [0,T]\times\Omega,
\end{align}
where $u$ is a scalar quantity, 
$a$ is a constant velocity field,
$\kappa>0$ is the diffusive constant, 
and $\Omega\subset \R$ is the one-dimensional domain with periodic boundary conditions.

The domain $\Omega$ is partitioned into $\Nel$ non-overlapping elements, and we define the mesh $\Tmeshh := \cup_{i=1}^\Nel \Elem_i$ 
by a finite collection of the elements $\Elem_i$. 
We denote the boundary of element $\Elem$ by $\partial \Elem$.
We let $\pTmesh_h := \LRc{\partial I:\Elem \in \Tmesh_h}$ be the collection of the boundaries of all elements.
For two neighboring elements $\Elem^+$ and $\Elem^-$ that share an interior
interface $\e = \Elem^+ \cap \Elem^-$, we denote by $q^\pm$ the trace of their
solutions on $\e$. 
We define $\nm$ as the unit outward normal vector on
the boundary $\partial \Elem^-$ of element $\Elem^-$, and we define $\np = -\nm$ as the unit outward
normal of a neighboring element $\Elem^+$ on $\e$.
On the interior face $\e$, we define the mean/average operator $\average{\bf v}$, where $\bf v$ is
either a scalar or a vector quantity, as
$\averageM{{\bf v}}:=\LRp{{\bf v}^- + {\bf v}^+}/2$, 
and the jump operator $\jumpL{\bf v} := {\bf v}^+ \cdot \np + {\bf v}^- \cdot \nm $.
Let $\Poly^{p}\LRp{D}$ denote the space of polynomials of degree at
most $p$ on a domain $D$. Next, we introduce discontinuous
piecewise polynomial spaces for scalars and vectors as
\begin{align*}
\Vhp{p}\LRp{\Tmeshh} &:= \LRc{v \in L^2\LRp{\Tmeshh}:
  \eval{v}_{\Elem} \in \Poly^p\LRp{\Elem}, \forall \Elem \in \Tmeshh},
\end{align*}
and a similar space $\Vhp{p}\LRp{\Elem}$ by replacing $\Tmeshh$ with $\Elem$. We define $\LRp{\cdot,\cdot}_\Elem$ as the $L^2$-inner product on an
element $\Elem \in \R^d$ and $\LRa{\cdot,\cdot}_{\partial \Elem}$ as the
$L^2$-inner product on the element boundary $\partial \Elem \in
\R^{d-1}$. Here, $d$ is the dimension. 
We define associated norm as $\norm{ \cdot }_{DG}:= \LRp{ \sum_{\Elem \in \Tmesh_h} \norm{ \cdot }_{\Elem}^2 }^\half$, where $\norm{ \cdot }_{\Elem}=\LRp{\cdot,\cdot}_\Elem^\half$.

The discontinuous Galerkin (DG) weak formulation of \eqnref{pde-convdiff-eq} yields the following: Seek $q,u \in \Vhp{p}\LRp{\Elem}$ such that 
\begin{subequations}
  \eqnlab{gov-convdiff-weak}
  \begin{align}
   \eqnlab{gov-convdiff-gradu}
   \LRp{\kappa^{-1}q,r}_\Elem &= 
     -\LRp{\dd{u}{x},r}_\Elem + \LRa{\nb(u - u^{**}),r}_{\partial \Elem}, \\
   \eqnlab{gov-convdiff-u}
   \LRp{\dd{u}{t},v}_\Elem
     &= -\LRp{ \dd{au+q}{x},v}_\Elem 
     + \LRa{ \nb \LRp{ au+q - au^{*} - q^{**} },v }_{\partial \Elem},
  \end{align}
 \end{subequations}
for all $r,v \in \Vhp{p}\LRp{\Elem}$ and each element $\Elem\in \Tmeshh$. 
Here we use the central flux for $u^{**}$ and $q^{**}$ for the diffusion operator, in other words, 
 $u^{**} = \averageM{\u}$  and $q^{**} = \averageM{\q}$, and the
Lax--Friedrich flux for the convection operator, $au^* = \averageM{au} + \half|a|\jumpL{u}$. We apply the nodal DG method \cite{hesthaven2007nodal} for discretizing \eqnref{gov-convdiff-weak}. The local solution is expanded by a linear combination of Lagrange interpolating polynomials. The Legendre--Gauss--Lobatto (LGL) points are used for interpolation. For integration in \eqnref{gov-convdiff-weak}, Gauss quadrature with $2(p+1)$ points is employed.

With explicit time integrators, 
we can eliminate the gradient $q$ term in \eqnref{gov-convdiff-gradu} by plugging it into \eqnref{gov-convdiff-u} 
and writing \eqnref{gov-convdiff-weak} as 
\begin{align}
  \eqnlab{gov-convdiff-ode} 
  \DD{u}{t} &= R (u),
\end{align}
where $R(u)$ is the spatial discretization of \eqnref{gov-convdiff-u}.
We let $u^{H} \in \Vhp{H}\LRp{\Elem}$ and $u^{L} \in \Vhp{L}\LRp{\Elem}$ 
respectively be the high- and the low-order approximations ($H>L$) on $\Elem$. 
The high- and the low-order spatial operators are represented by  $R^H(u^H)$ and $R^L(u^L)$, respectively. 

We introduce a low-order spatial filter $G$ such that $Gu^H \in \Vhp{L}(\Elem)$ 
by projecting the high-order approximation $u^H$ to the low-order approximation $u^L$ on each element, 
\begin{align}
\eqnlab{l2projection}
  \LRp{u^H - G u^H, v}_\Elem = 0, 
\end{align}
for all $v \in \Vhp{L}(\Elem)$. 
Applying the filter $G$ to the high-order semi-discretized form of \eqnref{gov-convdiff-ode} yields
\begin{align*}
  \DD{Gu^H}{t} &= GR^H (u^H). 
\end{align*}
As mentioned in \cite{de2022accelerating}, the tendency of $Gu^H$ is not the same as that of $u^L$; in other words, 
$\DD{Gu^H}{t} \ne \DD{u^L}{t}$. 
Our objective is to enhance the low-order solution $u^L$ by the filtered high-order solution $Gu^H$ 
through a neural network. To that end, we introduce the neural network source term in 
\eqnref{gov-convdiff-ode}, 
\begin{align}
\eqnlab{eq-convdiff-nnsource}
  \DD{\hat{u}^L}{t} &= R^L (\hat{u}^L) + S_\theta(\hat{u}^L). 
\end{align}
Here, $\hat{u}^L$  stands for 
the  prediction of the low-order operator with the neural network source term.

To obtain the training data, we generate multiple trajectories of the high-order solution $u^H$  on which the elementwise $L_2$ projection \eqnref{l2projection} is taken. In training, we randomly choose $n$-batch instances from the trajectories
and collect $m$ consecutive instances $\LRc{\LRp{Gu^H}_{t_s}^{(i)},\LRp{Gu^H}_{t_s+\dt}^{(i)},\cdots,\LRp{Gu^H}_{t_s+m\dt}^{(i)} }$ from each trajectory for $i=1,\cdots,n$.
Next, we predict the low-order approximation augmented by the neural network source term,  
$\LRc{\LRp{\hat{u}^L}_{t_s+\dt}^{(i)},\LRp{\hat{u}^L}_{t_s+2\dt}^{(i)},\cdots,\LRp{\hat{u}^L}_{t_s+m\dt}^{(i)} }$, 
by solving \eqnref{eq-convdiff-nnsource} with 
for $i=1,\cdots,n$; and we update the network parameters $\theta$ using the mean squared loss function,
\begin{align}
  \mathrm{loss} =  \frac{1}{n m}\sum_{i=1}^n \sum_{\ell=1}^m \norm{\LRp{Gu^H}_{t_s+\ell \dt}^{(i)} - \LRp{\hat{u}^L}_{t_s+\ell \dt}^{(i)}  }^2.  
  \eqnlab{lossfunc}
\end{align}

\subsection{Viscous Burgers' Equation}

Our third example is the viscous Burgers' equation defined on the time and space interval $(t,x) \in [0,T]\times\Omega$:
\begin{align}
  \eqnlab{pde-vburgers-eq}
 \dd{u}{t} + \half\dd{u^2}{x} 
 = \kappa \dd{u}{x^2} \text{ in } [0,T]\times\Omega,
\end{align}
where $u$ is a velocity field, $\kappa>0$ is the diffusive constant, and $\Omega\subset \R$ is the one-dimensional domain with periodic boundary conditions. 
The DG weak formulation is the same as \eqnref{gov-convdiff-weak} except for the convective term $\frac{u^2}{2}$ instead of $au$. 
The Lax--Friedrich flux for the convection operator is defined by $\LRp{\half u^2}^* = \averageM{\half u^2} + \half \tau \jumpL{u}$, where $\tau=\max\LRp{|u^+|,|u^-|}$. 

To obtain the training data, we generate a single trajectory of the high-order solution $u^H$ and take the elementwise $L_2$ projection \eqnref{l2projection} on $u^H$. In training, we randomly choose $n$-batch instances from the trajectory
and collect $m$ consecutive instances $\LRc{\LRp{Gu^H}_{t_s}^{(i)},\LRp{Gu^H}_{t_s+\dt}^{(i)},\cdots,\LRp{Gu^H}_{t_s+m\dt}^{(i)} }$ for $i=1,\cdots,n$. Then we update the network parameters $\theta$ using the same loss function in \eqnref{lossfunc}.

\section{Numerical Results}
\seclab{NumericalResults}
 In this section we illustrate our methodology for learning continuous source dynamics through numerical studies. 
 Our implementation leverages the automatic differentiation Python libraries   
 \texttt{PyTorch} \cite{paszke2019pytorch} and \texttt{JAX} \cite{jax2018github}. For the convection-diffusion equation and the viscous Burgers' equation, 
 the \texttt{Diffrax} \cite{kidger2022neural} NODE package and the \texttt{Optax} optimization library are used for training the neural network source term. 
We train the neural network source term on ThetaGPU at the Argonne Leadership Computing Facility using an NVIDIA DGX A100 node.
We measure the wall-clock times of the prediction of $\hat{u}^L$, $u^L$, and $u^H$ on Macbook Pro (i9 2.3 GHz).

\subsection{Two-Scale Lorenz 96 Model}

The training data is generated by using the ERK4 scheme and \eqnref{eq-lorenz96} with a timestep size of $\dt=0.005$. The simulation is warmed up for $t=3$ by using a randomly generated initial condition  $X\in \R^K$ and $Y\in \R^{JK}$.
After spinup, we set $t=0$ and integrate 
for $t \in [0,10]$. 
This process is repeated $300$ times. We take this generated data to be the ground truth solution.
  
  Next, we train the neural network source term in \eqnref{eq-lorenz96-slow-nnsource} using the training dataset, which is obtained by choosing random starting points in the truth dataset and integrating for five timesteps with a timestep size of $\dt$. 
We use the Adam optimizer \cite{kingma2014adam} with $2,000$ epochs and $100$ batches along with the loss function defined in \eqnref{lorenz96-loss}.

We plot the trained neural network source term $S_\theta$ and true coupling term $S$ in Figure \figref{L96-ss-source}. The neural network source term differs from the true coupling source term by less than $5\%$. In comparison with the true coupling term, the neural network source generally exhibits smoother profiles, as expected. 

Next we forecast the slow variables using the learned neural network source. 
We randomly choose an initial profile 
and integrate the model \eqnref{eq-lorenz96-slow-nnsource}
with $10\dt$.
\footnote{
The prediction timestep size $10\dt$ was ten times larger than the training timestep size $\dt$.
}
The true value $X$ (top) and the prediction $\hat{X}$ (middle) for the slow variables are shown in the Hovm\"{o}ller diagram in Figure \figref{L96-ss}.
The difference between the two simulations of $X$ and $\Xhat$ is shown in the bottom panel.
The prediction with the neural network source term is close to the true values up to $t=2$, which exceeds the reported Lyapunov time of the slow variables of the uncoupled model ($h=0$) that is $0.72$ \cite{bocquet2020bayesian}. For this specific initial condition, the augmenting neural network source model's predictability is about $2$ time units.

\begin{figure}[h!t!b!]
  \centering
  \includegraphics[trim=0cm 0cm 0cm 0cm,clip=true,width=0.9\textwidth]{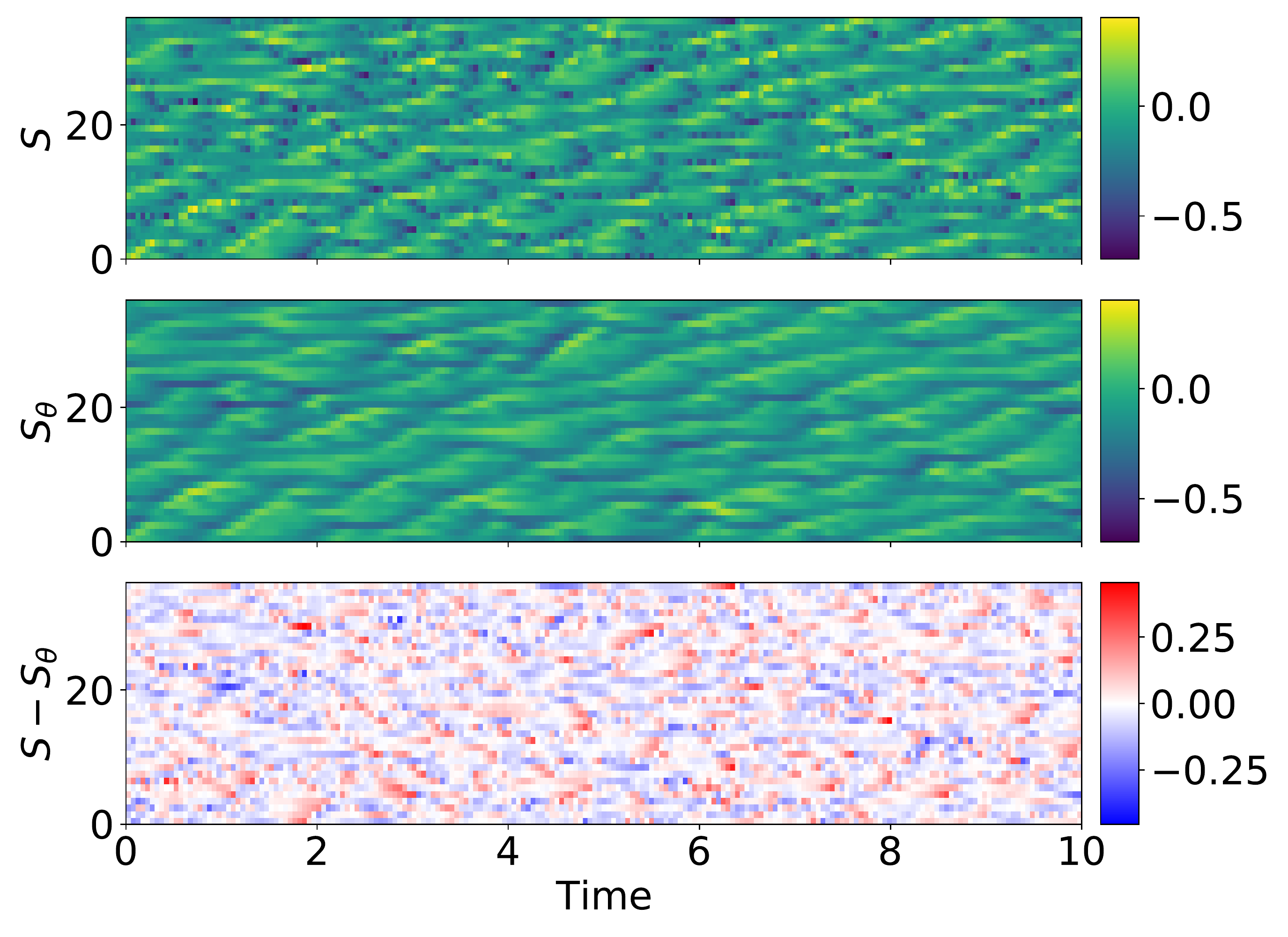}
  \caption{True source $S$ (top) and the trained source $S_\theta$ (middle)  shown for the two-scale Lorenz 96 model. 
  The bottom panel shows the difference between $S$ and $S_\theta$. 
  The trained source $S_\theta (X)$ is generated by using true $X$. 
  }
  \figlab{L96-ss-source}
\end{figure}

\begin{figure}[h!t!b!]
  \centering
  \includegraphics[trim=0cm 0cm 0cm 0cm,clip=true,width=0.9\textwidth]{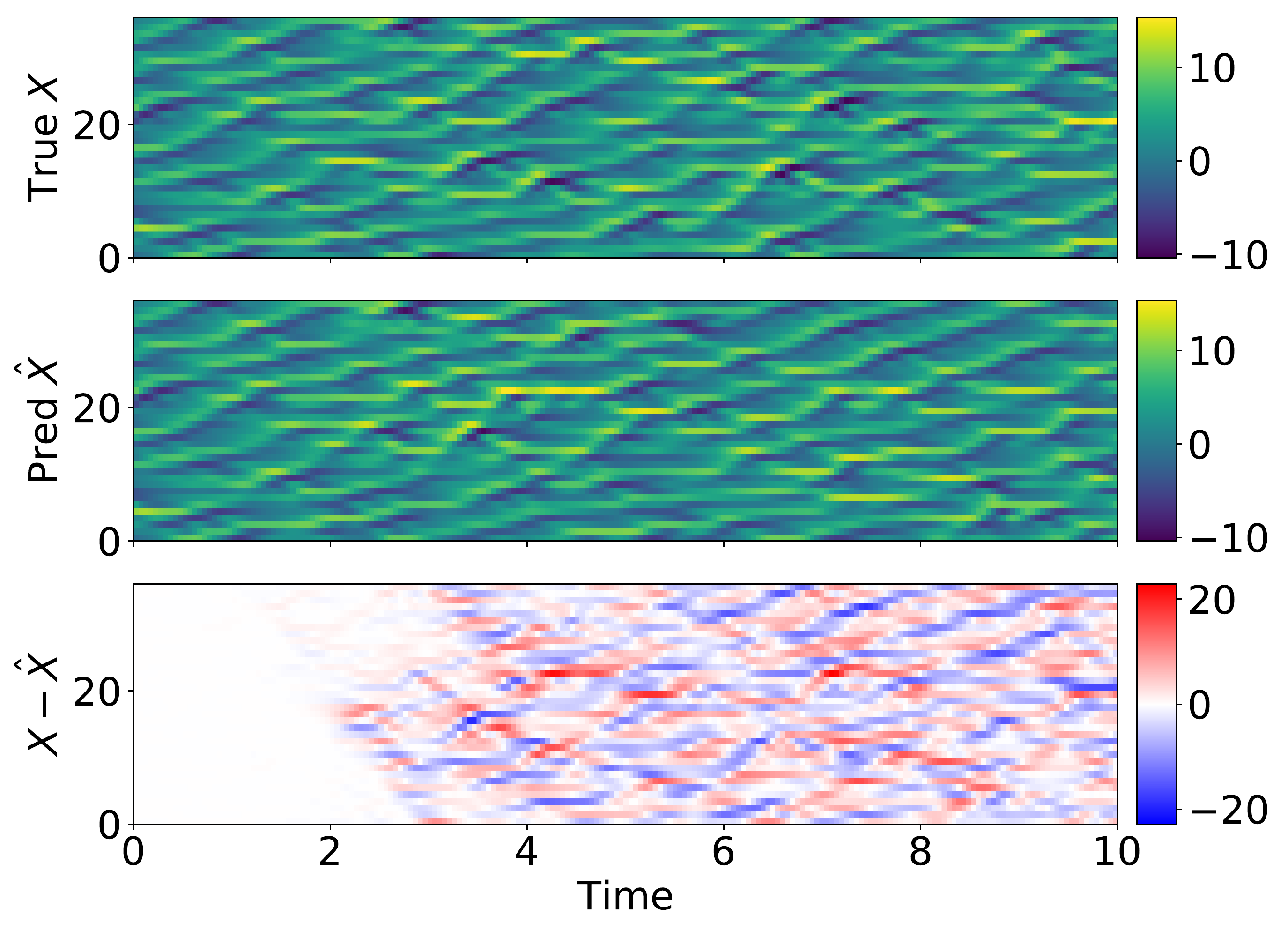}
  \caption{Hovm\"{o}ller diagram of the true model (top) and the trained model (middle) for the slow variables. 
  The bottom panel shows the difference between the two simulations.
  }
  \figlab{L96-ss}
\end{figure}

We also plot the Hovm\"{o}ller diagram with another initial condition in Figure 
\figref{L96-ss-xdt}.
We advance the trained model  \eqnref{eq-lorenz96-slow-nnsource} with the different timestep sizes of
$\LRc{1,2,5,10}\times\dt$. All the predictions show good agreement up to around $t=3$. 
As anticipated, the proposed method is not sensitive to the timestep sizes used for prediction because of the continuous aspect of the trajectory obtained by using a neural network source operator.  
In this example, the neural network source term enables the prediction of the slow variables with a step size that is $10$ times larger than the true model \eqnref{eq-lorenz96-slow}, without interacting with the fast variables.  
  
\begin{figure}[h!t!b!]
  \centering
  \includegraphics[trim=1.2cm 0.1cm 0.2cm 0.2cm,clip=true,width=0.9\textwidth]{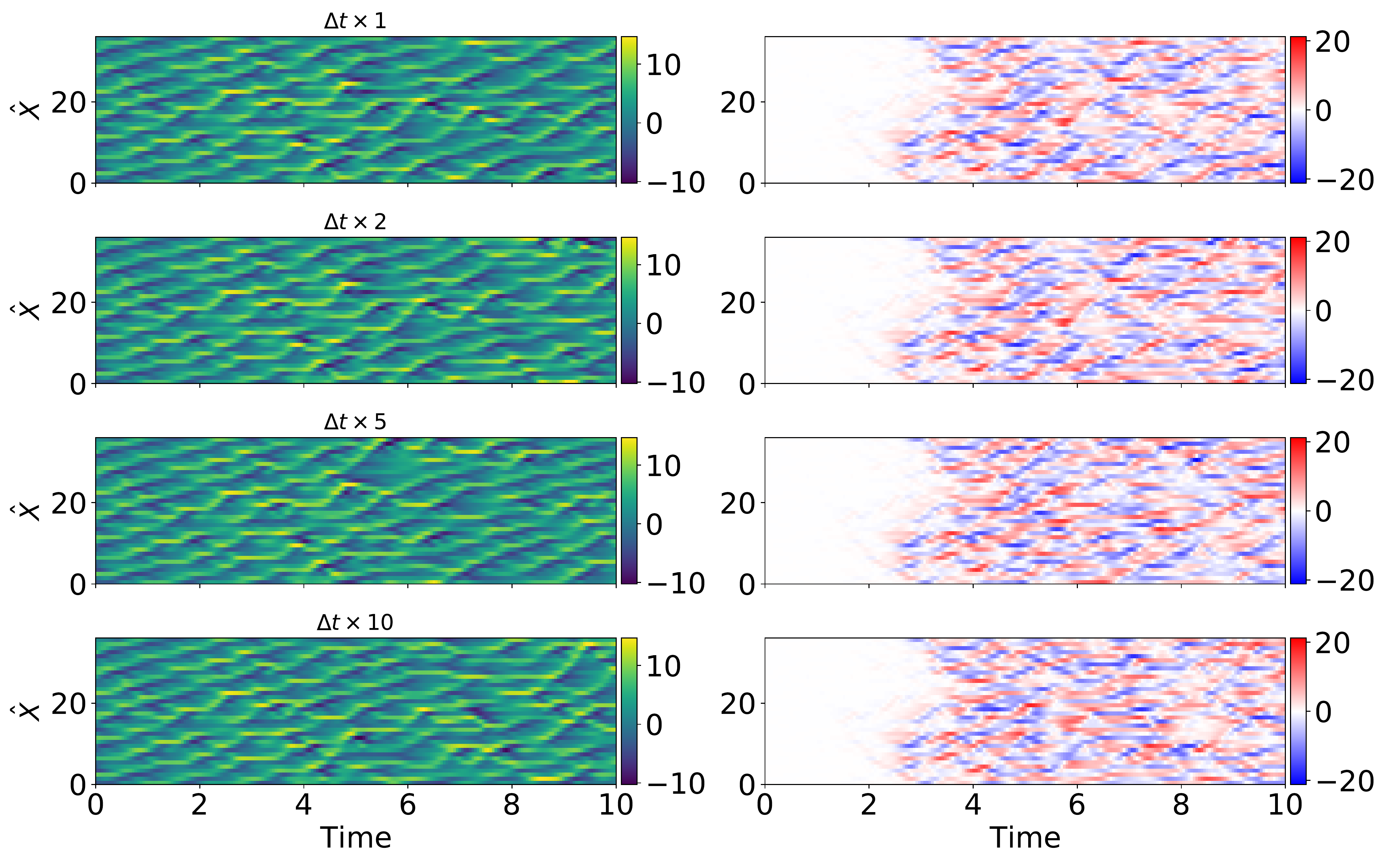}
  \caption{Hovm\"{o}ller diagram of the trained model (left) and the differences (right) with the true values for the slow variables. 
  Slow variables are advanced with the timestep sizes of $\LRc{1,2,5,10}\dt$.
  }
  \figlab{L96-ss-xdt}
\end{figure}

\subsection{One-Dimensional Convection-Diffusion Equation}

We next consider the propagation of an initial signal with different spectral components through the convection-diffusion equation on $\Omega=[0,1]$. The initial condition is therefore chosen as 
\begin{align}
\eqnlab{cd-ic}
u(t=0) =\sum_{i=1}^{4} \sin(2\pi \alpha_i (x-\phi)), 
\end{align}
with $\alpha_i=\{20,4,6,7\}$ and the phase function $\phi$ uniformly distributed from $0$ to $1$. We take $a=1$ and $\kappa=10^{-4}$ in \eqnref{pde-convdiff-eq}.



First, we run simulations using the ERK4 scheme with 100 randomly chosen phases $\phi$ and the timestep size of $10^{-4}$ for $t=[0,1]$ over the mesh of the fifth-order polynomial ($p=5$) and $50$ elements. 
Then, we project the time series data onto the first-order solutions $Gu^H$ by applying the filter $G$ on each element. For example, Figure \figref{cd-ic} shows the high-order solution $u^H$ (grey dashdot line) and its filtered low-order profile $Gu^H$ (brown solid line) at $t=0$. Next, 
we split the time series of $Gu^H$ into the training data for $t=[0,0.75]$ and 
the test data for $t=[0.75,1]$. 

For training the neural network source term, we randomly select $100$ batch instances of $Gu^H$ and integrate \eqnref{eq-convdiff-nnsource} for $5$ timesteps by using an explicit fifth-order RK method 
\cite{tsitouras2011runge} over the mesh with the first-order $p=1$ polynomial and $50$ elements and with $\dt=10^{-3}$, which is 10 times larger than the timestep that was used to compute the reference solution.
Since the first-order discontinuous Galerkin methods have two nodal points on each element, the neural network input and output layers have $100$ degrees of freedom.
We train the neural network source term by using the 
AdaBelief \cite{zhuang2020adabelief} optimizer with a learning rate of $10^{-4}$ and $3,000$ epochs. Figure \figref{cd-smallkappa-losshistory} shows the loss-epoch diagram. Both the training and test loss saturate at around $10^{-5}$.

\begin{figure}[h!t!b!]
  \centering
  \includegraphics[trim=0.5cm 0.5cm 0.1cm 0.8cm,clip=true,width=0.7\textwidth]{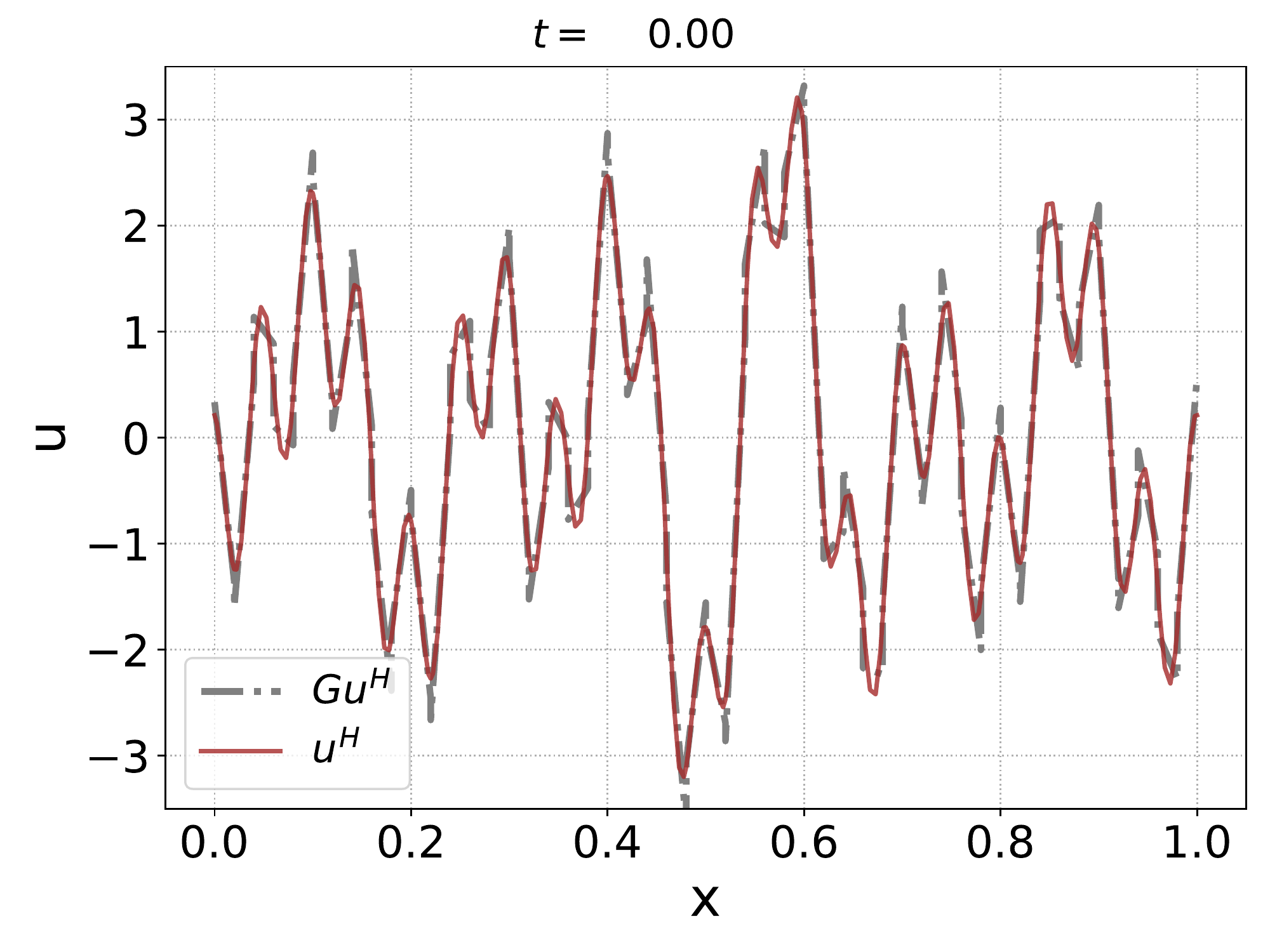}
  \caption{Convection-diffusion equation: the filtered (grey dashdot line) and the unfiltered (brown solid line) initial conditions. The fifth-order approximation of the initial condition is projected to the first-order polynomial approximation on each element. 
  }
  \figlab{cd-ic}
\end{figure}

\begin{figure}[h!t!b!]
  \centering
  \includegraphics[trim=0cm 0cm 0cm 0cm,clip=true,width=0.9\textwidth]{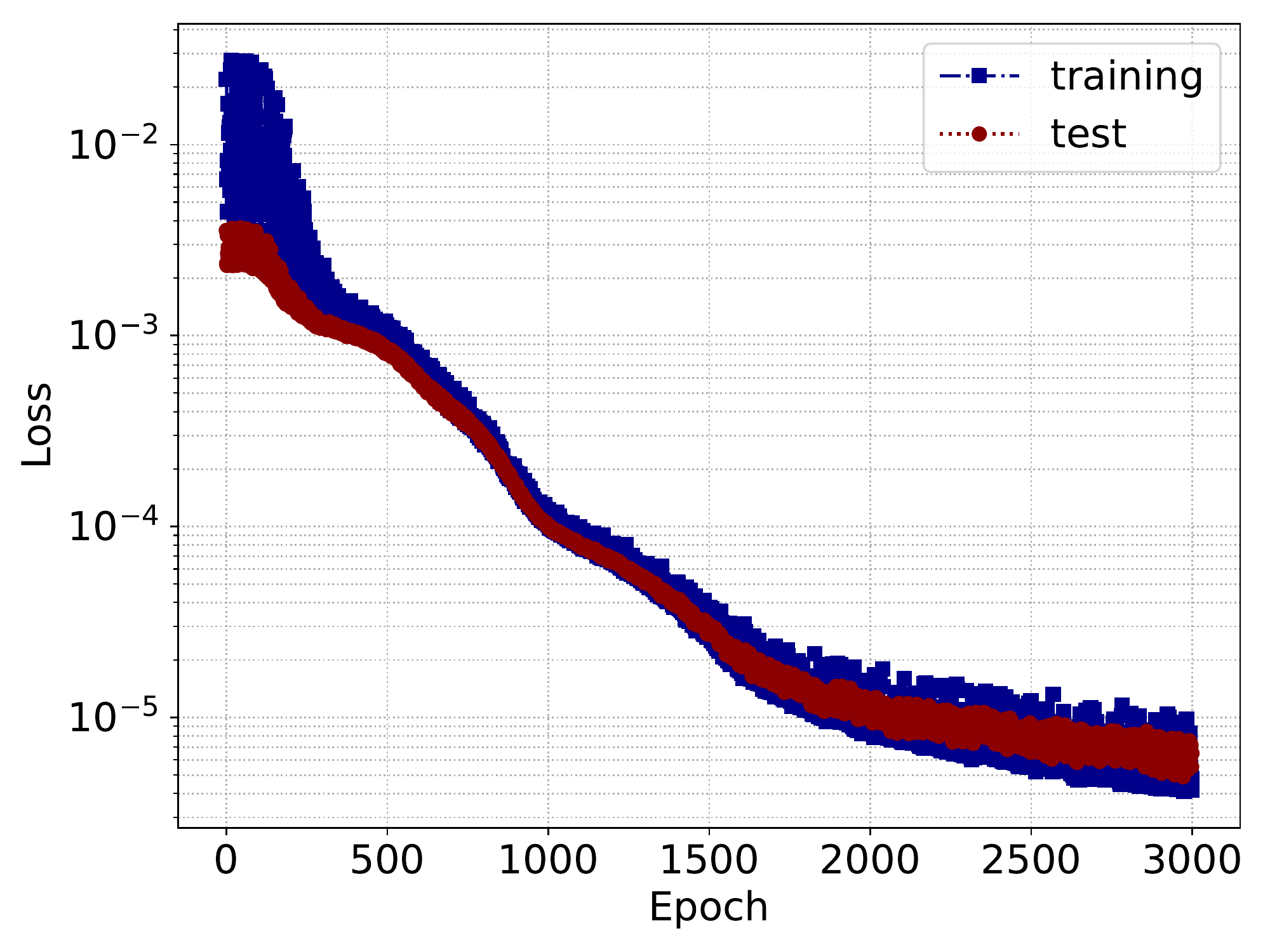}
  \caption{Convection-diffusion: the loss history computed with the training and test data.
  }
  \figlab{cd-smallkappa-losshistory}
\end{figure}

We let $u^L$ be the first-order DG solution without any neural network source term. We denote the first-order DG solution with the neural network source term by $\hat{u}^L$. 
In Figure \figref{cd-xtdiagram-smallkappa},  the $x$--$t$ diagrams depict how $u^L$ (top), $\hat{u}^L$ (middle), and $Gu^H$ (bottom) evolve in time. In comparison with the first-order DG solution $u^L$, 
the prediction with the neural network source $\hat{u}^L$ represents more accurately the filtered solution $Gu^H$. We compare $u^L$ and $G u^H$ as well as  $\hat{u}^L$ and $G u^H$, and we plot their differences in Figure \figref{cd-xtdiagram-smallkappa-diff}. We restrict the range of the colormap between $-1$ and $1$. The maximum absolute value of $u^L - Gu^H$ is $1.15$, whereas its $\hat{u}^L - Gu^H$ counterpart is $0.05$. The maximum $L_2$ errors of $u^L$ and $\hat{u}^L$ with respect to $Gu^H$ are $\norm{u^L - Gu^H}_{DG}=0.54$ and $\norm{\hat{u}^L - Gu^H}_{DG}=0.02$. 
The prediction with the neural network source term is therefore 25 times closer to the filtered solution $Gu^H$ than the low-order approximation $u^L$. 
Indeed, by adding the neural network source term, we note an enhanced accuracy of the first-order DG approximation.
In Figure \figref{cd-ss-smallkappa-t1} we show a snapshot of $u^L$, $\hat{u}^L$, and $Gu^H$ at the final time ($t=1$). We note that the neural network source term helps maintain a less smooth (polygonal) waveform of the filtered solution, when compared with the first-order DG method without the source term.
 
\begin{figure}[h!t!b!]
  \centering
  \includegraphics[trim=0.3cm 0.1cm 0.2cm 0.2cm,clip=true,width=0.9\textwidth]{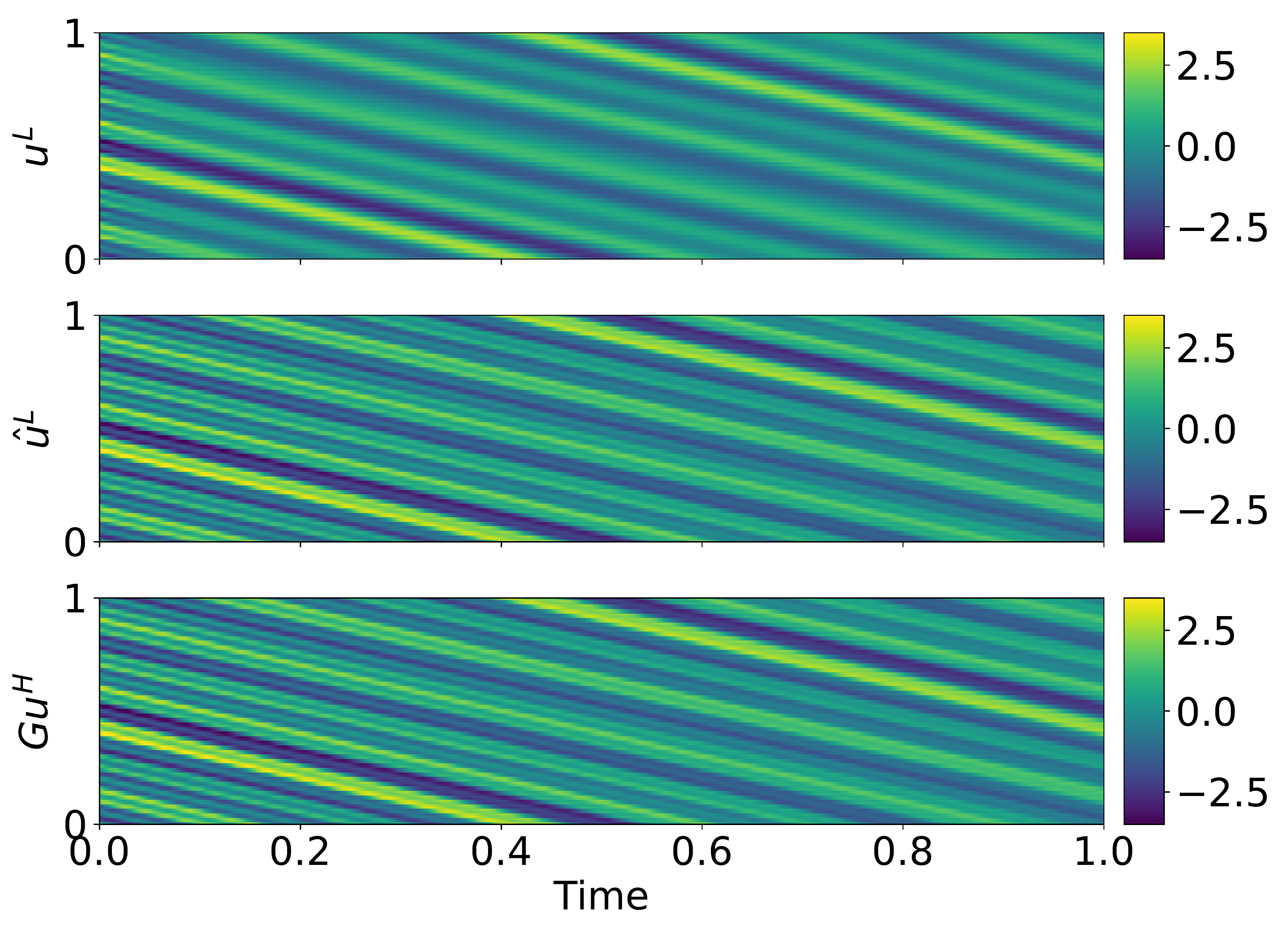}
  \caption{$x-t$ diagrams of the evolution of $u^L$ (top) , $\hat{u}^L$ (middle), and $Gu^H$ (bottom) for $\kappa=10^{-4}$. 
  }
  \figlab{cd-xtdiagram-smallkappa}
\end{figure}

\begin{figure}[h!t!b!]
  \centering
  \includegraphics[trim=0.3cm 0.5cm 0.2cm 0.2cm,clip=true,width=0.9\textwidth]{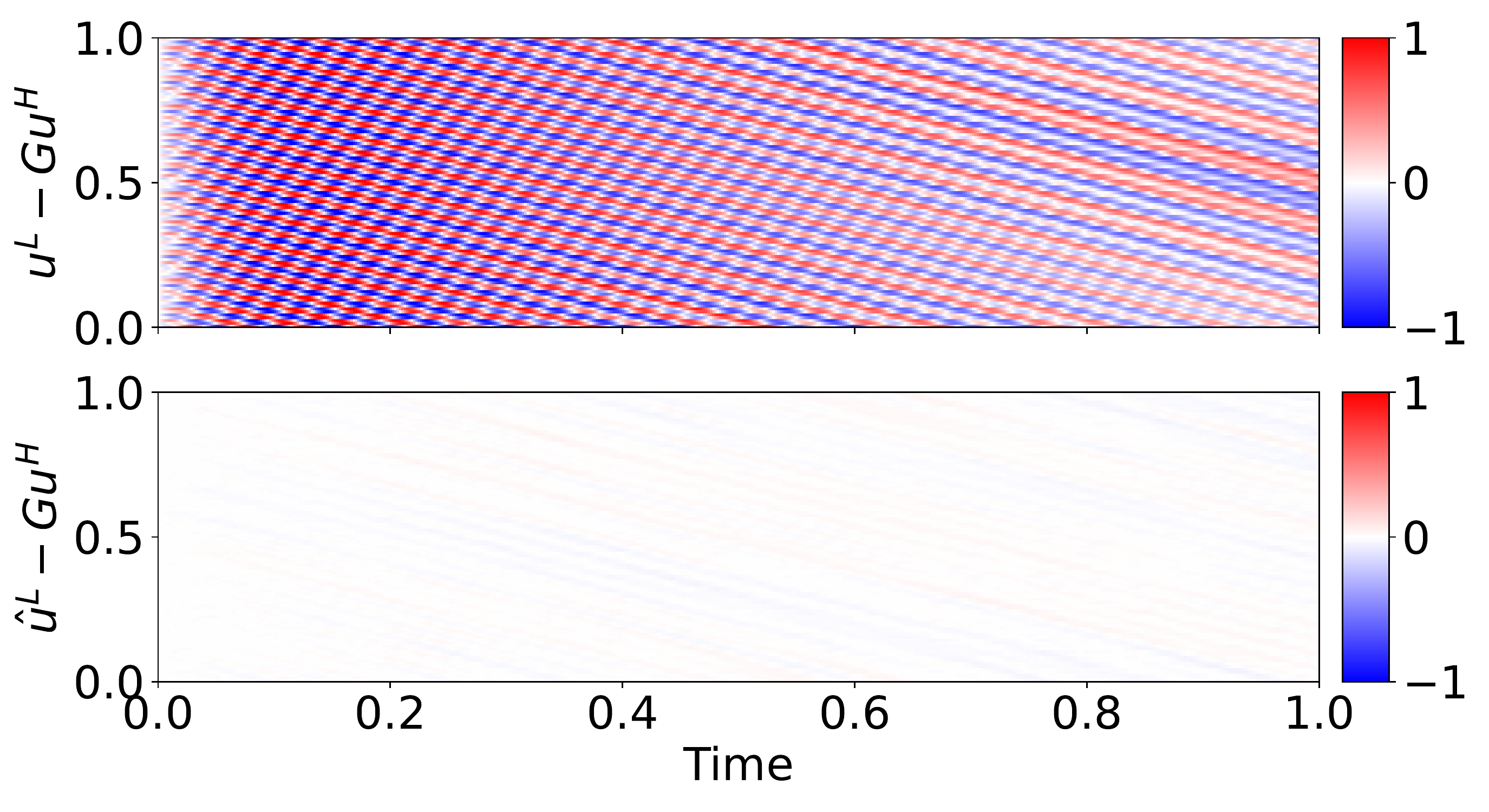}
  \caption{Difference between $u^L$ and $G u^H$ (top)
  and the difference between $\hat{u}^L$ and $G u^{H}$ (bottom) for $\kappa=10^{-4}$.
  }
  \figlab{cd-xtdiagram-smallkappa-diff}
\end{figure}

\begin{figure}[h!t!b!]
  \centering
  \includegraphics[trim=0.3cm 0.5cm 0.2cm 0.2cm,clip=true,width=0.9\textwidth]{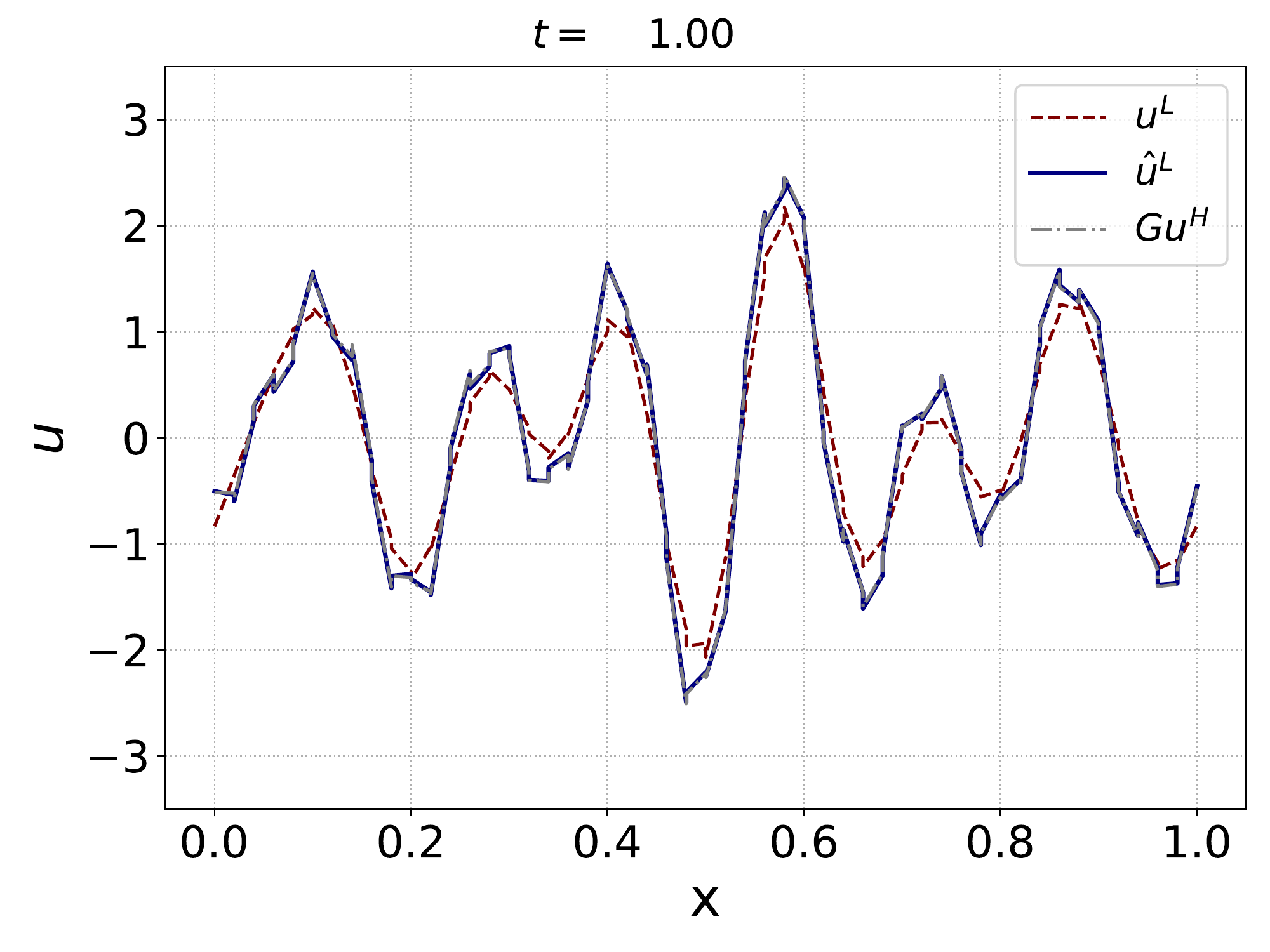}
  \caption{Snapshots of $u^L$ (top), $\hat{u}^L$ (middle), and $G u^H$ (bottom) at $t=1$ for $\kappa=10^{-4}$.
  }
  \figlab{cd-ss-smallkappa-t1}
\end{figure}

Now we consider the case with $\kappa=10^{-3}$. 
To account for the new diffusive coefficient, we create further training and test data.
Because of the increasing viscosity, the initial condition \eqnref{cd-ic} quickly diffuses as time passes. Throughout the evolution, the relative gap between $u^L$ and $G u^H$ rapidly decreases, as illustrated in Figure \figref{cd-xtdiagram-diff}.
At the final time ($t=1$), $u^L$ looks similar to $\hat{u}^L$, as shown in Figure \figref{cd-ss-t1}.
However, the difference between $\hat{u}^L$ and $G u^H$ appears to be smaller than that of $u^L$ and $G u^H$. This is illustrated in Figure \figref{cd-xtdiagram-diff}, where we show a range between $-0.5$ and $0.5$. The maximum absolute value of $u^L - Gu^H$ is $0.55$, and the 
 $\hat{u}^L - Gu^H$ counterpart is $0.06$. The maximum $L_2$ errors of $u^L$ and $\hat{u}^L$ with respect to $Gu^H$ are $\norm{u^L - Gu^H}_{DG}=0.23$ and $\norm{\hat{u}^L - Gu^H}_{DG}=0.02$.
We observe that $\hat{u}^L$ is approximately ten times closer to $Gu^H$ than $u^L$.

\begin{figure}[h!t!b!]
  \centering
  \includegraphics[trim=0.3cm 0.5cm 0.2cm 0.2cm,clip=true,width=0.9\textwidth]{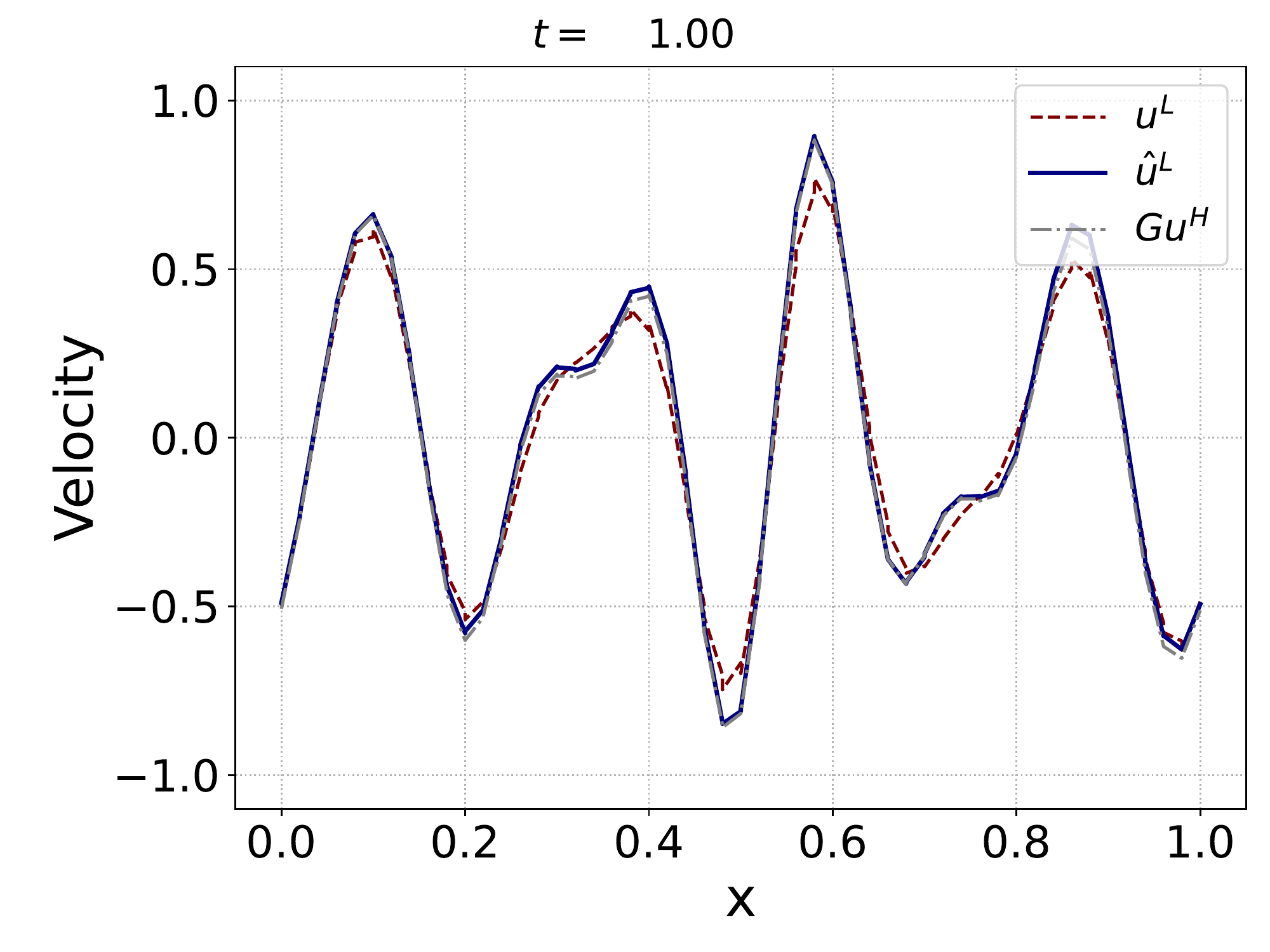}
  \caption{Snapshots of $u^L$ (top), $\hat{u}^L$ (middle), and $G u^H$ (bottom) at $t=1$ for $\kappa=10^{-3}$.
  }
  \figlab{cd-ss-t1}
\end{figure}

\begin{figure}[h!t!b!]
  \centering
  \includegraphics[trim=0.3cm 0.1cm 0.2cm 0.2cm,clip=true,width=0.9\textwidth]{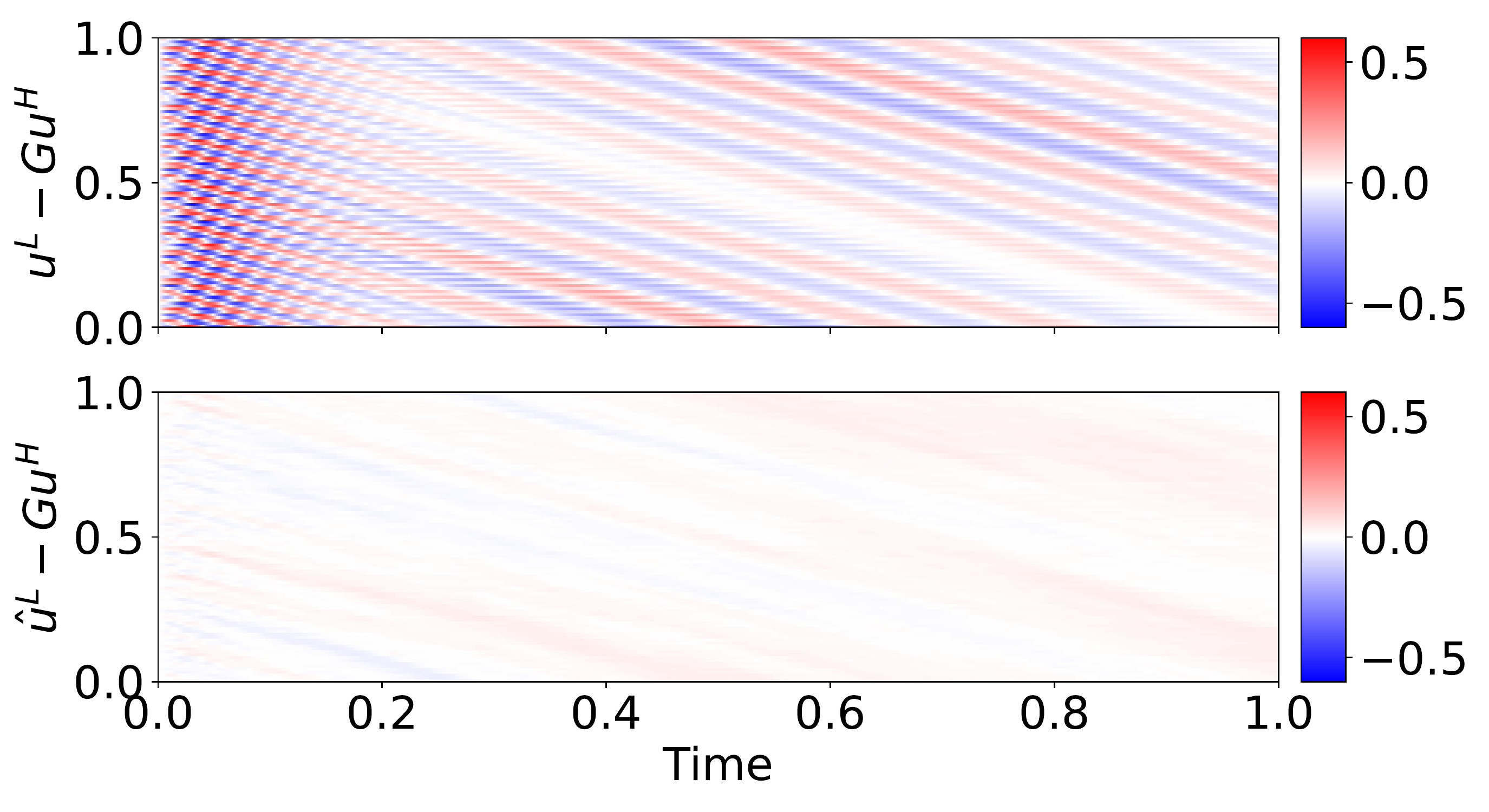}
  \caption{Difference between $u^L$ and $G u^H$ (top)
  and the difference between $\hat{u}^L$ and $G u^{H}$ (bottom) for $\kappa=10^{-3}$.
  }
  \figlab{cd-xtdiagram-diff}
\end{figure}


Next we find the polynomial equivalent to the augmented approximation. 
We denote the second- and the third-order DG approximations by $u_2$ and $u_3$, respectively. We let $u^L_2$ and $u^L_3$ respectively be the projected solutions of $u_2$ and $u_3$ onto the first-order polynomial space on each element. 
For the initial conditions of $u_2$ and $u_3$, we interpolate the initial profile of $u^L$ to the second- and the third-order polynomial spaces. We integrate $u^L$, $\hat{u}^L$, $u_2$, and $u_3$ with $\dt=10^{-3}$ for $t=[0,1]$ over
the mesh of $50$ elements. 
We measure the relative $L_2$ errors of $u^L$, $\hat{u}^L$, $u_2^L$, and $u^L_3$ by taking $Gu^H$ as the ground truth in Figure \figref{cd-relerr-history}. For $\kappa=10^{-4}$, the augmented DG approximation $\hat{u}^L$ is comparable to the projected third-order DG approximation $u^L_3$. For $\kappa=10^{-3}$, because of the increased diffusive coefficient, the gap between $u^L_2$ and $Gu^H$ quickly shrinks.
The enhanced DG approximation $\hat{u}^L$ is not comparable to the projected second-order DG approximation $u^L_2$. However,  
it still exhibits reduced relative error
in comparison with the first-order DG approximation.


\begin{figure}[h!t!b!]
  \centering
  
\subfigure[$\kappa=10^{-4}$]{  \includegraphics[trim=0.0cm 0.5cm 0.2cm 0.2cm,clip=true,width=0.43\textwidth]{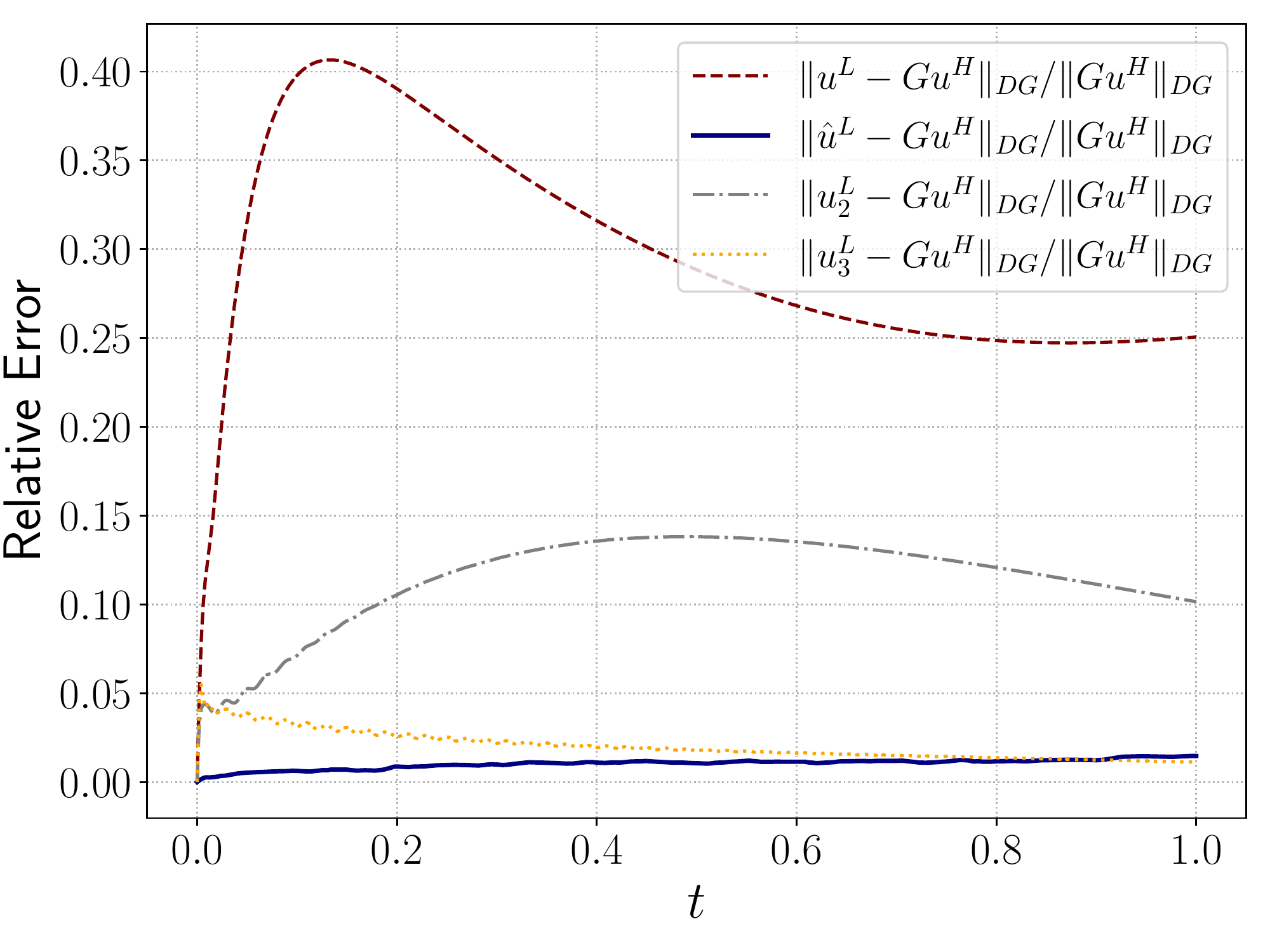}
}
\subfigure[$\kappa=10^{-3}$]{  \includegraphics[trim=0.0cm 0.5cm 0.2cm 0.2cm,clip=true,width=0.43\textwidth]{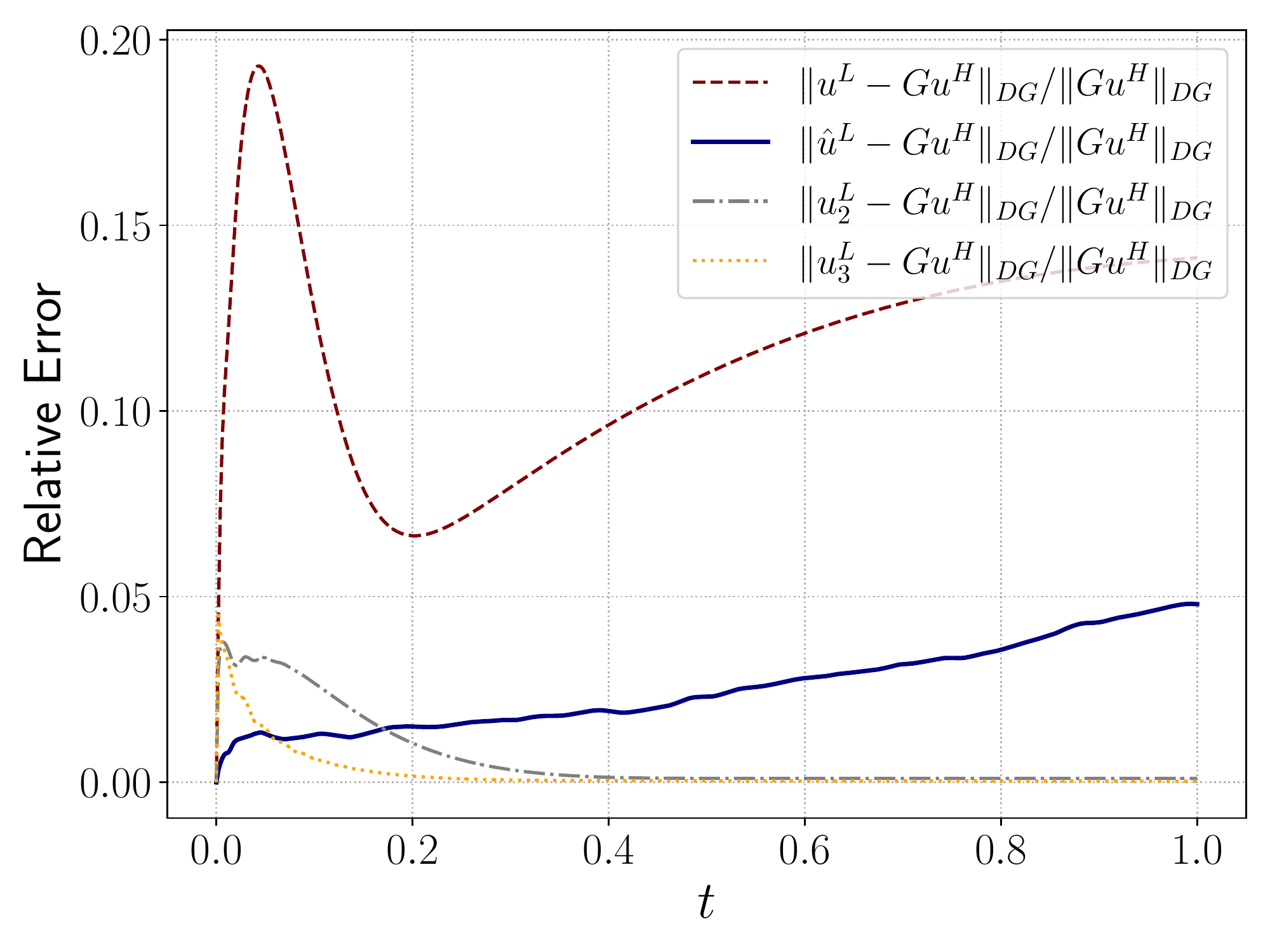}
}

  \caption{History of relative errors of $u^L$, $\hat{u}^L$, $u^L_2$, and $u^L_3$ for (a) $\kappa=10^{-4}$ and (b) $\kappa=10^{-3}$. The filtered solution $G u^H$ is considered as a reference solution.
  }
  \figlab{cd-relerr-history}
\end{figure}

For the predictions of $u^H$, $u^L$, $\hat{u}^L$, $u^L_2$, and $u^L_3$, we additionally provide the wall-clock time in Table \tabref{cd-walltime}. 
Because of the just-in-time (JIT) compilation property of \texttt{JAX}, 
we measure the wall clock two times. 
The reason is that 
when a function is first called in \texttt{JAX}, it will get compiled, and the resulting code will get cached. As a result, the second run proceeds much faster than the first. We simulate the case with $\kappa=10^{-4}$ ten times for the same time interval $[0,1]$ and report the average wall-clock time for the first and the second run.
We choose the stable timestep sizes of $\dt^H=0.001$ and $\dt^L=0.009$, which correspond to the largest steps that do not lead to a blowup. The corresponding Courant numbers for the convective term $Cr_a=\frac{a \dt}{\dx}$ are 1.36 and 1.44, respectively. The stable timestep sizes of $\dt =0.0045$ and $\dt=0.0025$ are taken for $u^L_2$ and $u^L_3$. 
Here, $\dx$ is the smallest distance between Legendre--Gauss--Lobatto points. 
For the first run with JIT compiling time, $\hat{u}^L$ takes slightly more time than the others, but $u^H$, $u^L$, $u^L_2$, and $u^L_3$ are comparable to one another.
For the second run with a cached code, 
the wall clock of the low-order solution $u^L$ is $12$ times smaller than that of the high-order solution $u^H$. The neural network augmented low-order solution $\hat{u}^L$ is about three times faster than the high-order solution $u^H$. The wall clock of $\hat{u}^L$ is four times larger than that of the low-order solution $u^L$ and comparable to that of the projected third-order DG approximation $u^L_3$. 

\begin{table}[t] 
    \caption{Wall-clock time for the predictions of $u^H$, $u^L$, $\hat{u}^L$, $u^L_2$, and $u^L_3$ for the convection-diffusion equation with $\kappa=10^{-4}$. }
    \tablab{cd-walltime} 
    \begin{center} 
    \begin{tabular}{*{1}{c}|*{1}{c}|*{1}{c}|*{1}{c}} 
    \hline 
      & $\dt $  & JIT Compile wc~[$ms$] & Simulation wc~[$ms$] \tabularnewline 
    \hline\hline 
    $u^H$ &   $0.001 $       &    749 & 16.4 \tabularnewline
    $u^L$ &   $0.009$        &    787 & 1.3  \tabularnewline
    $\hat{u}^L$ & $0.009$    &    996 & 4.8  \tabularnewline
    $u^L_2$ &   $0.0045$     &    793 & 2.4  \tabularnewline
    $u^L_3$ &   $0.0025$     &    781 & 4.8  \tabularnewline
    \hline\hline 
    \end{tabular} 
    \end{center} 
  \end{table}

\subsubsection*{Comparison with the discrete corrective forcing approach}

Our proposed method can be seen as a generalization of the discrete corrective forcing approach introduced in \cite{de2022accelerating} because in our approach we learn the underlying continuous source dynamics through the NODE instead of a post-correction.  In this section we briefly review the discrete corrective forcing approach and compare it with our proposed one. 

The discrete corrective forcing approach has two steps: advance the solution with the coarse-grid operator $R^L$ and then correct the low-order approximation $u^L$ using a corrective forcing in a postprocessing manner,  
\begin{subequations}
\begin{align}
  U_{n,i}^L &= u_n^L + \dt^L\sum_{j=1}^{i-1} a_{ij} R^L_j, \quad i=1,2,\cdots,s,\\
  u_{n+1}^L &= u_n^L + \dt^L \sum_{i=1}^s b_i R^L_i,\\
  \eqnlab{discrete-corrective-forcing}
  \tilde{u}_{n+1} &= u_{n+1}^L + \dt^L S_\theta(u_n^L),
\end{align}
\end{subequations}
where $\dt^L$ is the coarse timestep size and $\tilde{u}$ is the corrected approximation. The correction in \eqnref{discrete-corrective-forcing} can be viewed as a split timestepping step equivalent to a forward Euler method, which has first-order accuracy in time.  

We let $\dt^H$ be the timestep size used for generating training and test data of $Gu^H$. We take $\dt^L=10\times \dt^H$ and consider $\kappa=10^{-4}$. 
According to the procedure in \cite{de2022accelerating}, 
we train the corrective forcing term as follows. First we  generate 
$S_n = \frac{G u^H_{n+1} - u^L_{n+1}}{\dt^L} $ at $t_n^L$. 
Then we construct  a nonlinear map from the input feature of $Gu^H_n$ to the output feature of $S_n$ by using a feed-forward neural network. The discrete corrective forcing term is again trained by using the  
AdaBelief \cite{zhuang2020adabelief} optimizer with a learning rate of  $10^{-3}$, $100$ batches, and $3000$ epochs.
Note that we use the same neural network architecture for both the discrete and the continuous corrective forcing approaches. The work in \cite{de2022accelerating} considers the dependency of previous filtered solutions ($Gu^H_{n-r}$,$Gu^H_{n-r+1}$,$\cdots$,$Gu^H_{n}$), but for simplicity we do not consider the dependency of the earlier filtering solutions.

Next we predict both the $\tilde{u}^L$ and $\hat{u}^L$ for $t \in [0,1]$ with $\dt^L=10\dt^H$. We compare the $x$--$t$ diagrams of the two approaches and plot them in 
Figure  \figref{cd-xtdiagram-lala-10dt-diff}. Both the discrete and the continuous corrective forcing approaches show good agreement with each other as well as the training and test data of $Gu^H$. The maximum difference of $|\tilde{u} - Gu^H|$ is $0.03$, and the maximum difference of $|\hat{u} - Gu^H|$ is $0.05$. The discrete-in-time approach exhibits marginally higher accuracy than does the continuous-in-time approach.

\begin{figure}[h!t!b!]
  \centering
  \includegraphics[trim=0.3cm 0.1cm 0.2cm 0.2cm,clip=true,width=0.95\textwidth]{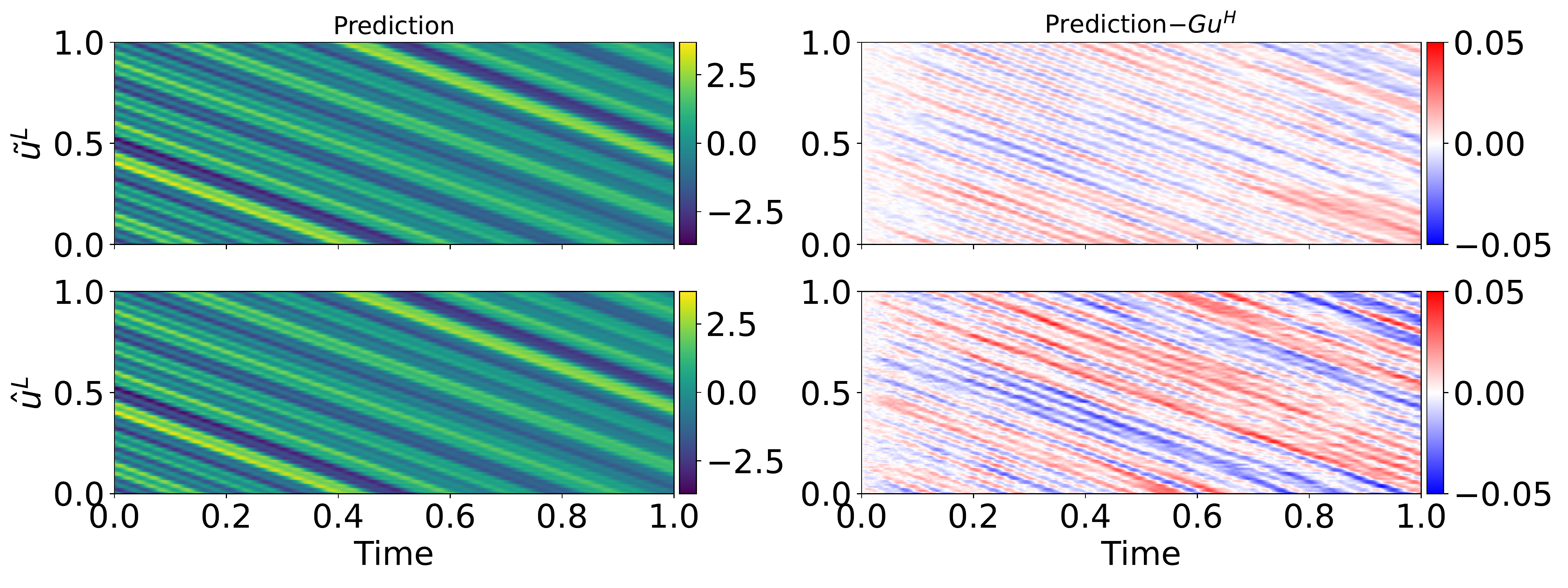}
  \caption{Difference between $\tilde{u}^L$ and $G u^H$ (top),
  and the difference between $\hat{u}^L$ and $G u^{H}$ (bottom).
  }
  \figlab{cd-xtdiagram-lala-10dt-diff}
\end{figure}

For the case with $\dt^L=20\times \dt^H$, however, we  observe that the continuous corrective forcing approach outperforms the discrete corrective forcing approach, as shown in Figure \figref{cd-xtdiagram-lala-20dt-diff}. 
The maximum difference of $|\tilde{u} - Gu^H|$ is about ten times larger than that of $|\hat{u} - Gu^H|$. 
\footnote{
$\max|\tilde{u} - Gu^H|=0.63$ and 
$\max|\hat{u} - Gu^H|=0.05$.
}
This is expected because the discrete corrective forcing approach learns a map from the low-order solution to the corrective forcing approximation with a specific timestep size and the first-order Euler method. Thus, predicting with different $\dt_L$ results in poor performance. 
On the other hand, the neural ODE approach learns the corrective forcing operator at a continuous level. This leads to consistent results that do not degrade the simulation predictability by changing the timestep size from $10\dt$ to $20\dt$. 

\begin{figure}[h!t!b!]
  \centering
  \includegraphics[trim=0.3cm 0.1cm 0.2cm 0.2cm,clip=true,width=0.95\textwidth]{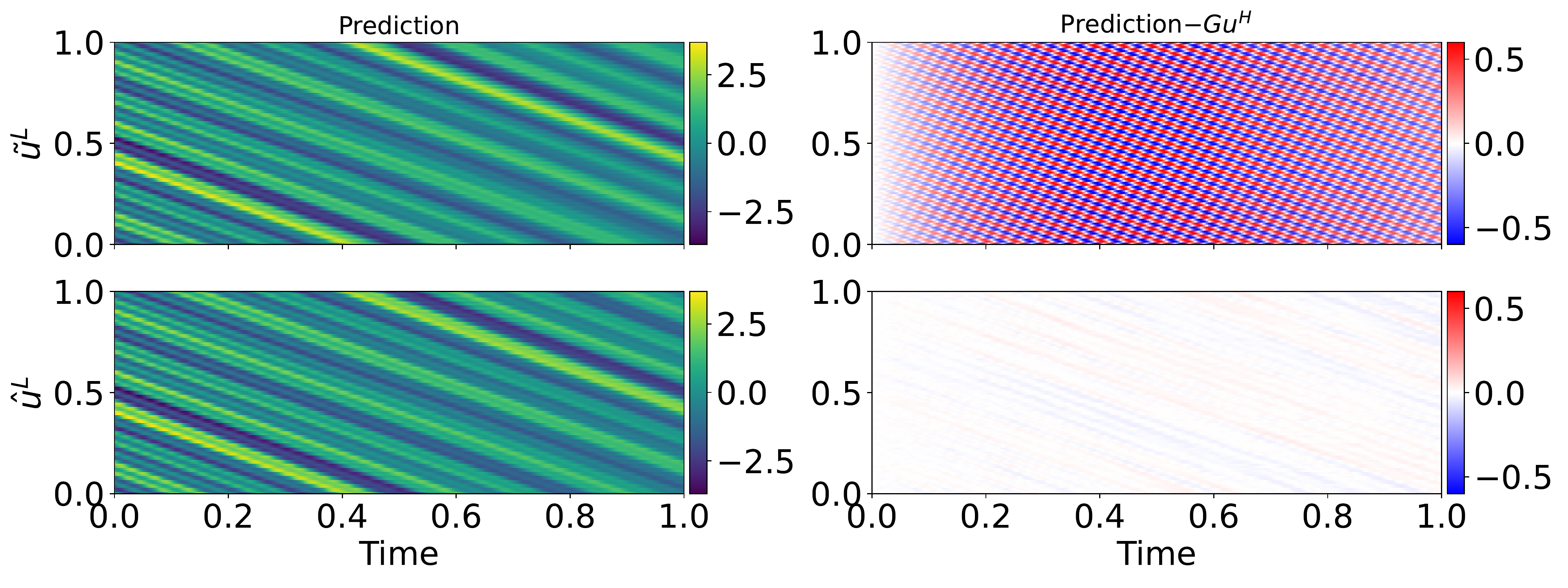}
  \caption{Difference between $\tilde{u}^L$ and $G u^H$ (top), 
  and the difference between $\hat{u}^L$ and $G u^{H}$ (bottom).
  }
  \figlab{cd-xtdiagram-lala-20dt-diff}
\end{figure}

To see the sensitivity of the prediction with respect to timestep size, we predict $\hat{u}^L$ and $\tilde{u}^L$
with timestep sizes of $\dt^L=\LRc{1,2,5,10,20,50}\times 10^{-4}$, and 
we measure the  relative $L_2$ errors at $t=0.5$ and $t=1$ by taking $Gu^H$ as the ground truth. 
Figure \figref{cd-relerr} shows the relative errors of $\frac{\norm{\hat{u}^L - Gu^H}_{DG}}{\norm{ Gu^H}_{DG}}$ (red dot line) and $\frac{\norm{\tilde{u}^L - Gu^H}_{DG}}{\norm{Gu^H}_{DG}}$  (blue dashdot line)  against timestep sizes. In general the continuous corrective forcing approach demonstrates  better accuracy than does the discrete corrective forcing approach. The error of $\hat{u}^L$ converges to a certain nonzero constant with decreasing $\dt^L$, for  example,  
$\frac{\norm{\hat{u}^L - Gu^H}_{DG}}{\norm{ Gu^H}_{DG}} \rightarrow 0.011$ at $t=0.5$ and
$\frac{\norm{\hat{u}^L - Gu^H}_{DG}}{\norm{ Gu^H}_{DG}} \rightarrow 0.015$ at $t=1$. 
The error of $\tilde{u}^L$, however, does not show any convergence; rather, it has the minimum relative error at $\dt^L=10^{-3}$ in which discrete corrective forcing was trained.

\begin{figure}[h!t!b!]
  \centering
   
\subfigure[$t=0.5$]{  \includegraphics[trim=0.0cm 0.5cm 0.2cm 0.2cm,clip=true,width=0.43\textwidth]{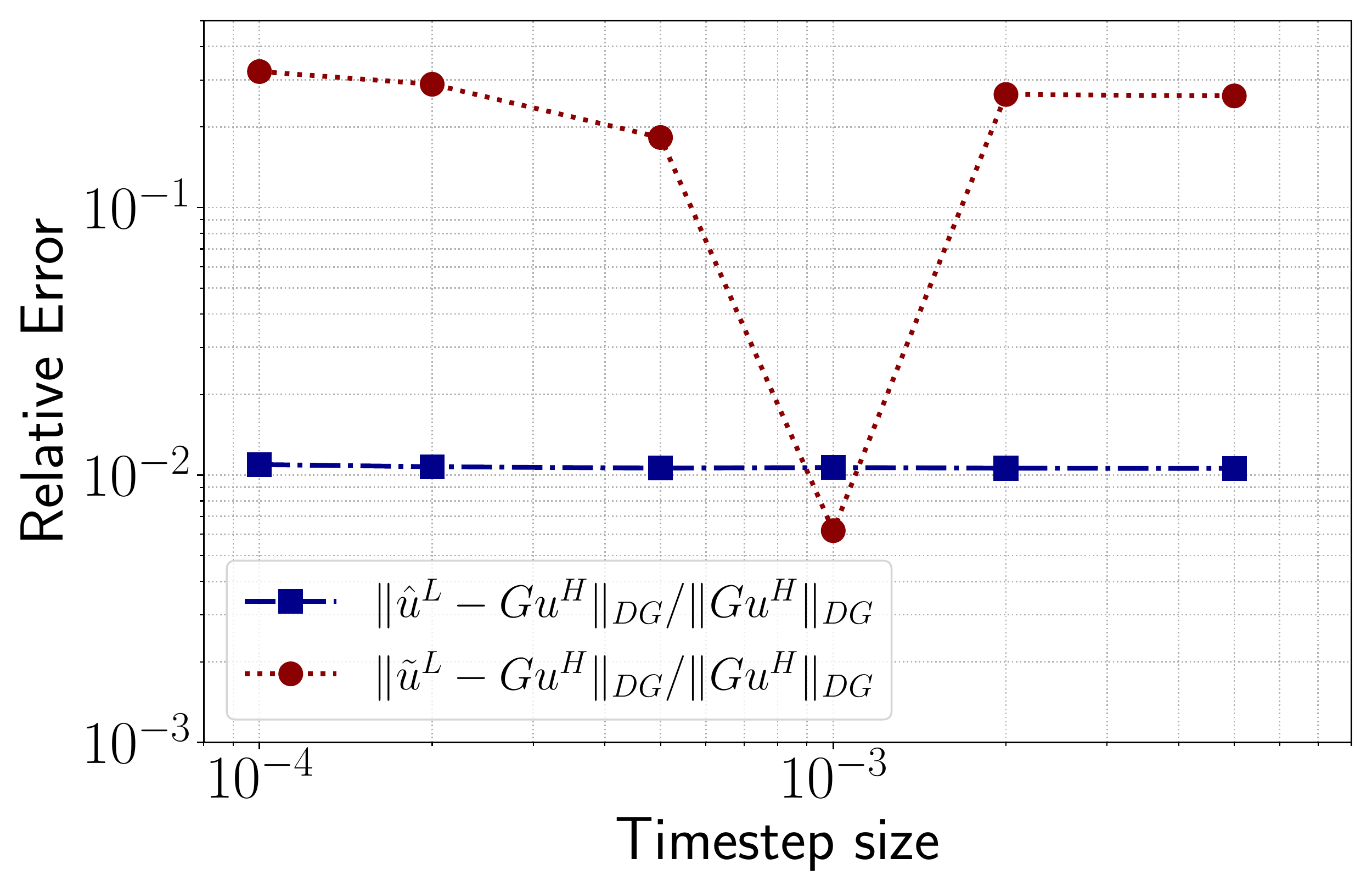}
}
\subfigure[$t=1.0$]{  \includegraphics[trim=0.0cm 0.5cm 0.2cm 0.2cm,clip=true,width=0.43\textwidth]{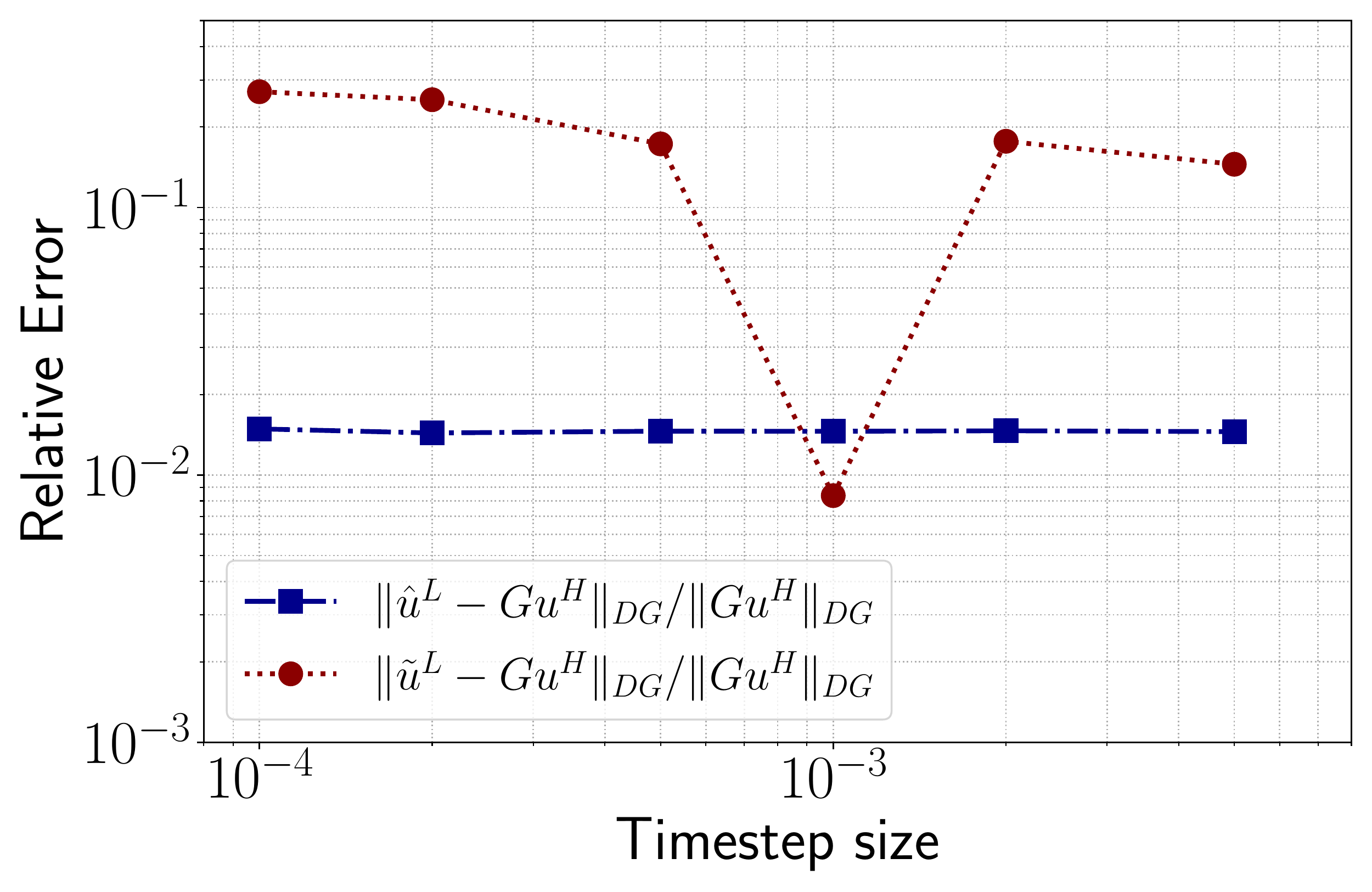}
}

  \caption{Relative errors of $\hat{u}^L$ and $\tilde{u}^L$ at (a) $t=0.5$ and (b) $t=1$. The filtered solution $G u^H$ is considered as a reference solution.
  }
  \figlab{cd-relerr}
\end{figure}

\subsection{One-Dimensional Viscous Burgers' Equation}





The kinetic energy is transferred from large-scale structures to smaller and smaller scales as a fluid moves. 
Energy cascade is a basic characteristic of turbulent flows that can be seen in Navier--Stokes equations. 
The viscous Burgers' equation is a simplified form of Navier--Stokes equations while keeping similar energy cascade to the Navier--Stokes system. 
Thus, the one-dimensional Burgers' equation is widely used for studying turbulence and sub-grid scale modeling \cite{love1980subgrid,basu2009can,frisch2002burgulence,labryer2015framework,san2018neural}.

We consider an example of Burgers turbulence \cite{bec2007burgers,maulik2018explicit} in the domain of $x\in[0,2\pi]$ with a periodic boundary condition. An initial energy spectrum is defined by 
$$ E(k) = A_0k^4 \exp\LRp{-\frac{k}{k_0}^2}$$
with $A_0=\frac{2k_0^{-5}}{3\sqrt{\pi}}$ and the peak wave number $k_0$ in the spectrum. The Fourier coefficients of high-order solution $u^H$ are obtained by 
$$ \mathfrak{u}(k) = N\sqrt{2 E(k)}\cos(2\pi \phi(k)),$$
where $k$ is the wave number in spectral domain, $N$ is the number of evenly spaced grid points, and $\phi$ is the phase function uniformly distributed from $0$ to $1$. 
The initial condition is generated by taking the inverse Fourier transform with $k_0=10$ and $N=32768$ and interpolating to the Legende--Gauss--Lobatto points, 
$u^H(t=0,x) = \Ical^H \LRp{\mathcal{F}^{-1}\LRc{\mathfrak{u}}}$. 
\footnote{
The interpolation is performed on each element by using the multiquadratic radial basis function $\phi(r)=\sqrt{1+(r/\epsilon)^2}$ in \texttt{scipy.interpolate.Rbf}.
}
The viscosity $\kappa$ is chosen as $0.005$. 
We integrate the viscous Burgers' model by using the ERK4 scheme with the timestep size of $\dt^H=5\times 10^{-4}$ and project the trajectory onto the first-order solutions $Gu^H$.  
Figure \figref{Burgers-burgulence-ic} shows the snapshots of the velocity field of $u^H$ at $t=[0,0.5,1]$ over the mesh of $\Nel=64$ elements and the eighth-order polynomial ($p=8$). 
\begin{figure}[h!t!b!]
  \centering
  \includegraphics[trim=0cm 0cm 0cm 0cm,clip=true,width=0.5\textwidth]{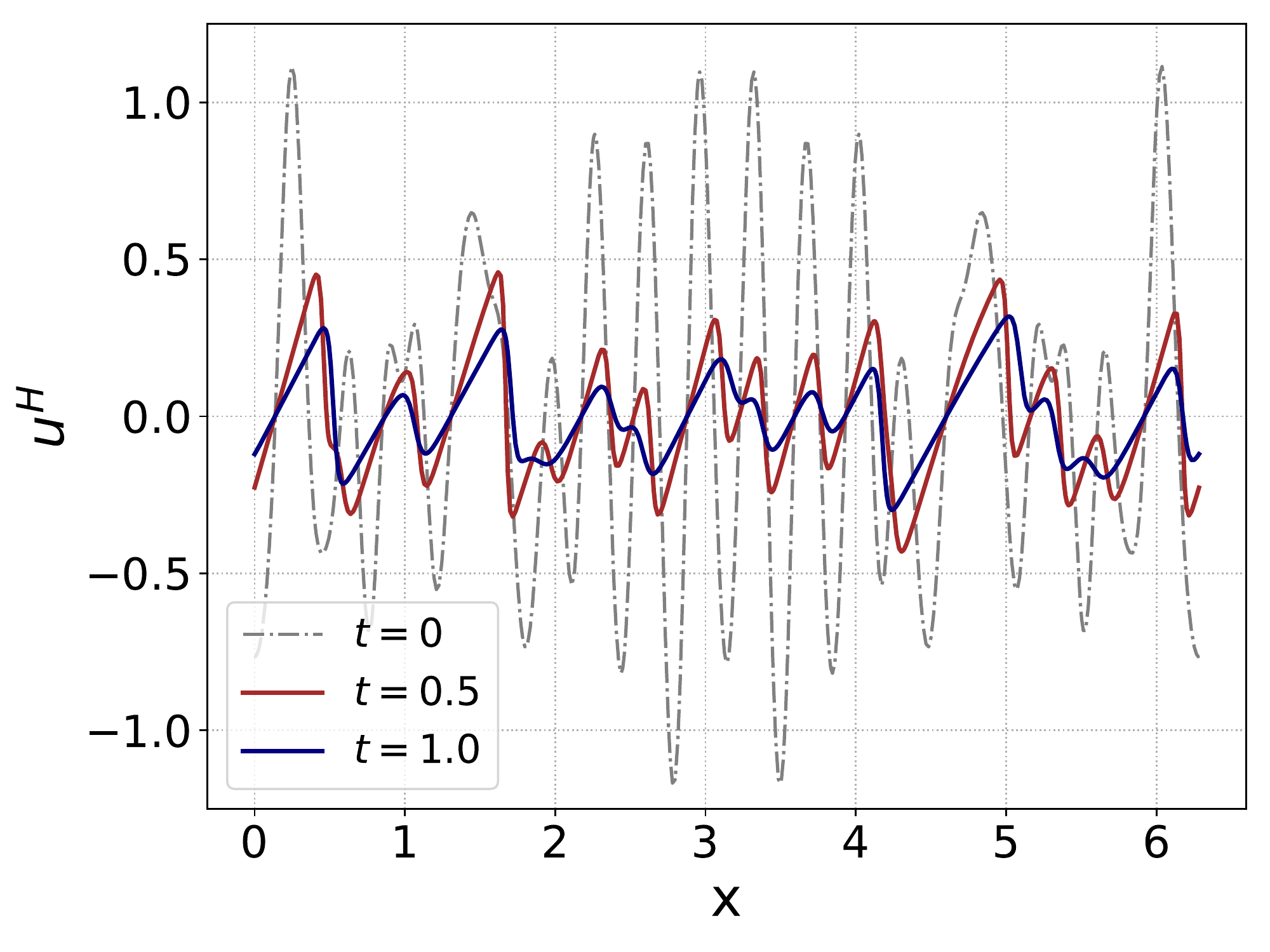}
  \caption{Burgulence: the snapshots of the velocity field of the high-order solution $u^H$ at $t=0$ and $t=1$. 
  }
  \figlab{Burgers-burgulence-ic}
\end{figure}
Next, 
we split the time series of $Gu^H$ into the training data for $t=[0,0.75]$ and 
the test data for $t=[0.75,1]$. 

For training the neural network source term, we randomly select $100$ batch instances of $Gu^H$ and integrate \eqnref{eq-convdiff-nnsource} for $5$ timesteps by using an explicit fifth-order RK method 
\cite{tsitouras2011runge} over the mesh with the first-order $p=1$ polynomial and $64$ elements and with $\dt=2.5\times 10^{-2}$, which is $50$ times larger than the timestep used for the filtered solution $Gu^H$. The neural network input and output layers have $128$ degrees of freedom.
We train the neural network source term by using the 
AdaBelief \cite{zhuang2020adabelief} optimizer with a learning rate of $10^{-3}$ and $500$ epochs, as shown in Figure \figref{cd-smallkappa-losshistory}. 

\begin{figure}[h!t!b!]
  \centering
  \includegraphics[trim=0cm 0cm 0cm 0cm,clip=true,width=0.6\textwidth]{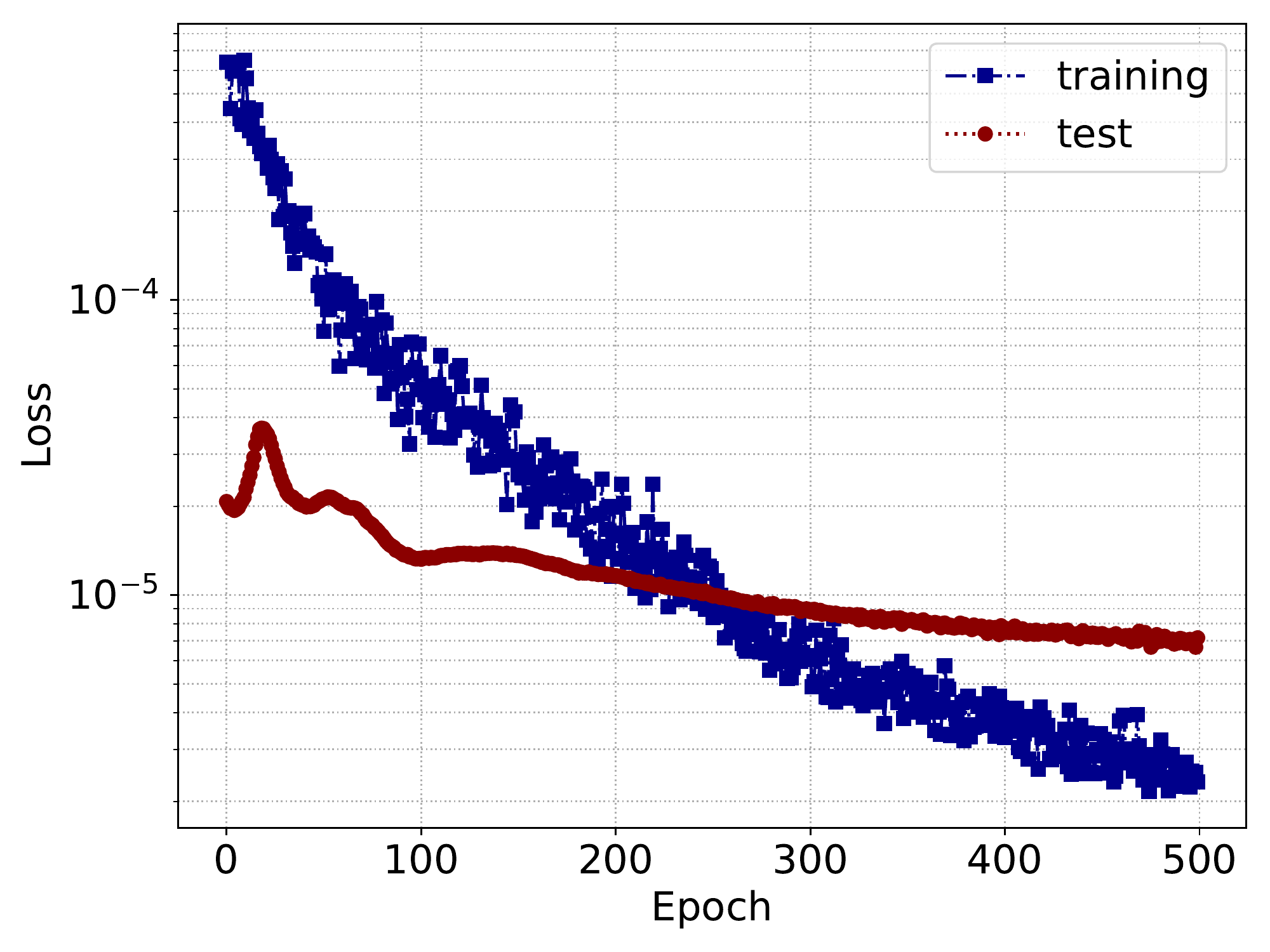}
  \caption{Burgulence: the loss history computed with the training and test data.
  }
  \figlab{Burgulence-losshistory}
\end{figure}

With the trained neural network, we ran a simulation over the mesh of $p=1$ and $\Nel=64$ with $\dt^L=50\times \dt^H$. 
Figure \figref{burgers-Burgulence-ss} shows the snapshot of $u^L$, $\hat{u}^L$, and $Gu^H$ at $t=0.5$ and $t=1.0$. 
We see that the neural network source term improves the low-order solution accuracy.
\begin{figure}[h!t!b!]
  \centering
\subfigure[$t=0.5$]{  \includegraphics[trim=0.0cm 0.0cm 0.0cm 0.8cm,clip=true,width=0.43\textwidth]{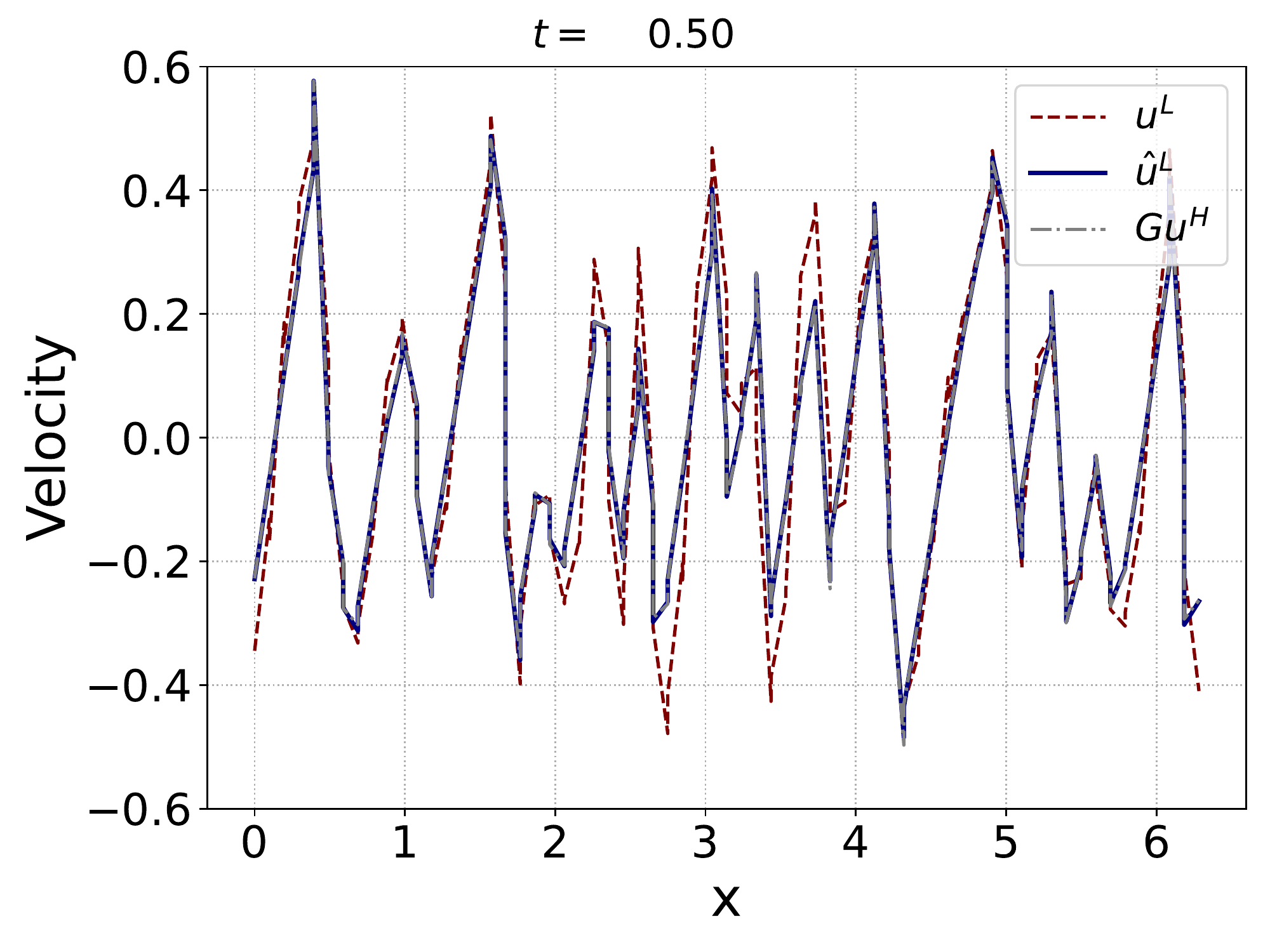}
}
\subfigure[$t=1.0$]{  \includegraphics[trim=0.0cm 0.0cm 0.0cm 0.8cm,clip=true,width=0.43\textwidth]{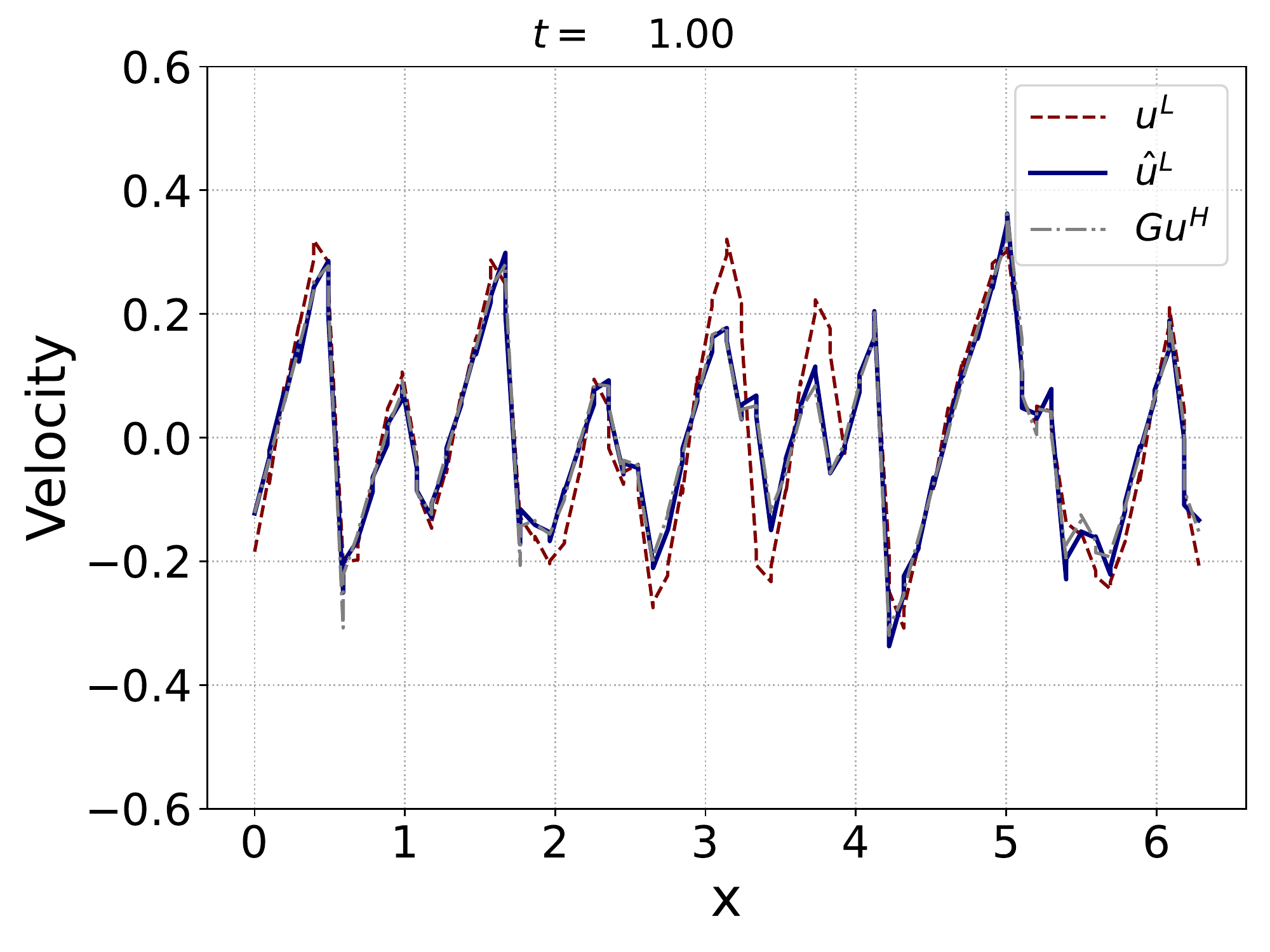}
}
  \caption{Burgulence: snapshots of $u^L$, $\hat{u}^L$, and $Gu^H$ at (a) $t=0.5$ and (b) $t=1.0$. 
  }
  \figlab{burgers-Burgulence-ss}
\end{figure}
We can also see that the augmenting neural network source term improves the accuracy of the kinetic energy spectrum in Figure \figref{burgers-spectra}. 
We first interpolate a DG solution on a uniform mesh $\Ical(u_{DG})$ with $N$ points, then we take the Fourier transform of that result, and finally we compute the angle averaged energy spectrum, 
$E(k)=\half\LRp{\hat{E}(k) + \hat{E}(-k)}$ for $k>0$.
Here, $\hat{E}(k)=\frac{1}{2N^2} \abs{\Fcal^{-1}\LRc{\Ical(u_{DG} )}}^2$. 
We take $N=64$ for the first-order DG solutions ($u^L$, $Gu^H$, $\hat{u}^L$), and $N=512$ for the eight-order ($u^H$) DG approximation. Two vertical lines with the cutoff wave numbers $k=32$ and $k=256$ are displayed in Figure \figref{burgers-spectra}. 
The grey solid line represents the kinetic energy spectrum of the unfiltered high-order solution $u^H$. Ideal scaling $k^{-2}$ \cite{maulik2018explicit} is represented by the blue dashed line. 
The green dotted line indicates the filtered high-order solution $Gu^H$.
The energy spectrum of $\hat{u}^L$ (navy solid line) is well overlapped with the spectrum of the filtered solution $Gu^H$, 
whereas the energy spectrum of $u^L$ (red dashdot line) departs slightly from it. 

\begin{figure}[h!t!b!]
  \centering
   
\subfigure[$t=0.5$]{  \includegraphics[trim=0.0cm 0.0cm 0.0cm 0.0cm,clip=true,width=0.43\textwidth]{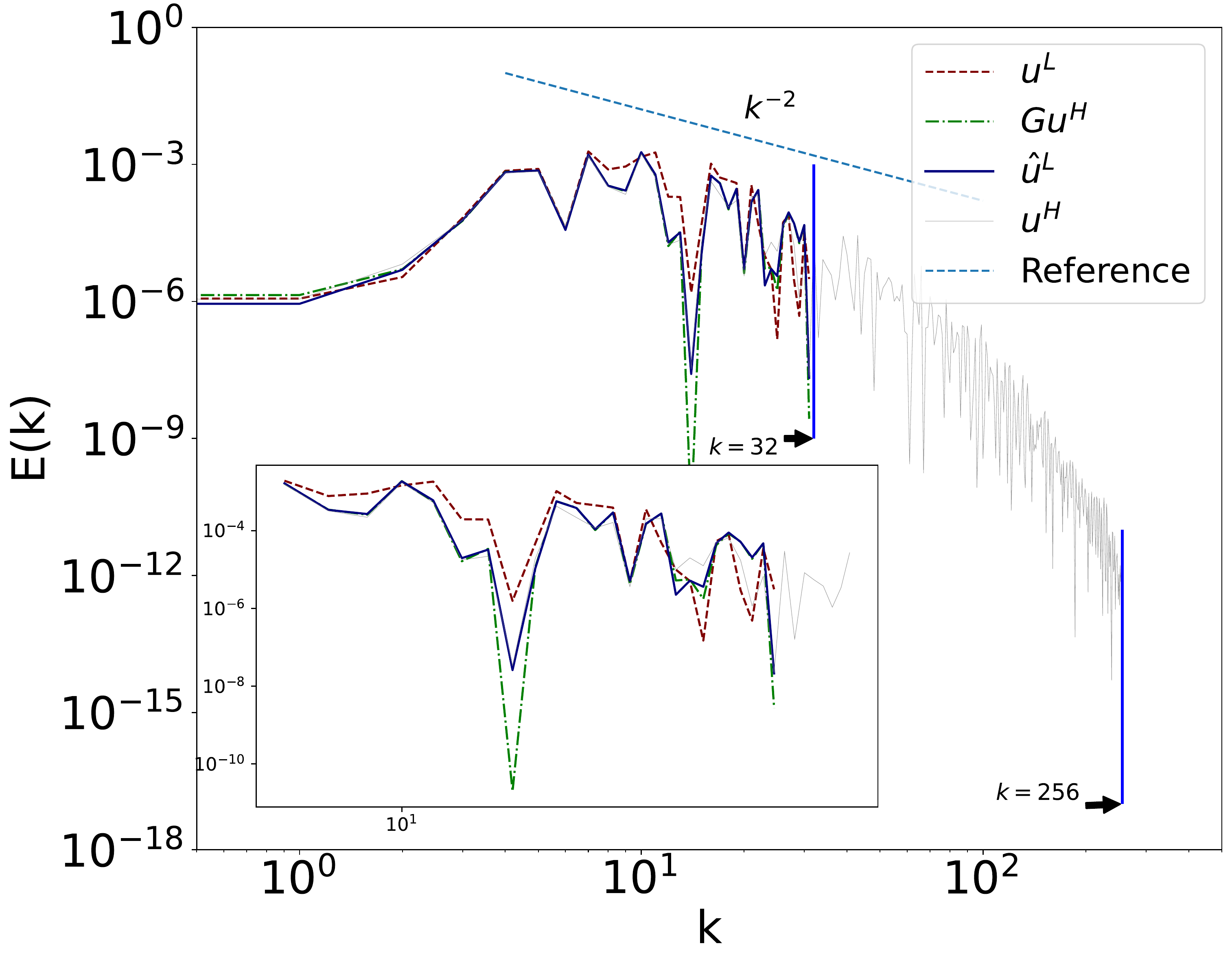}
}
\subfigure[$t=1.0$]{  \includegraphics[trim=0.0cm 0.0cm 0.0cm 0.0cm,clip=true,width=0.43\textwidth]{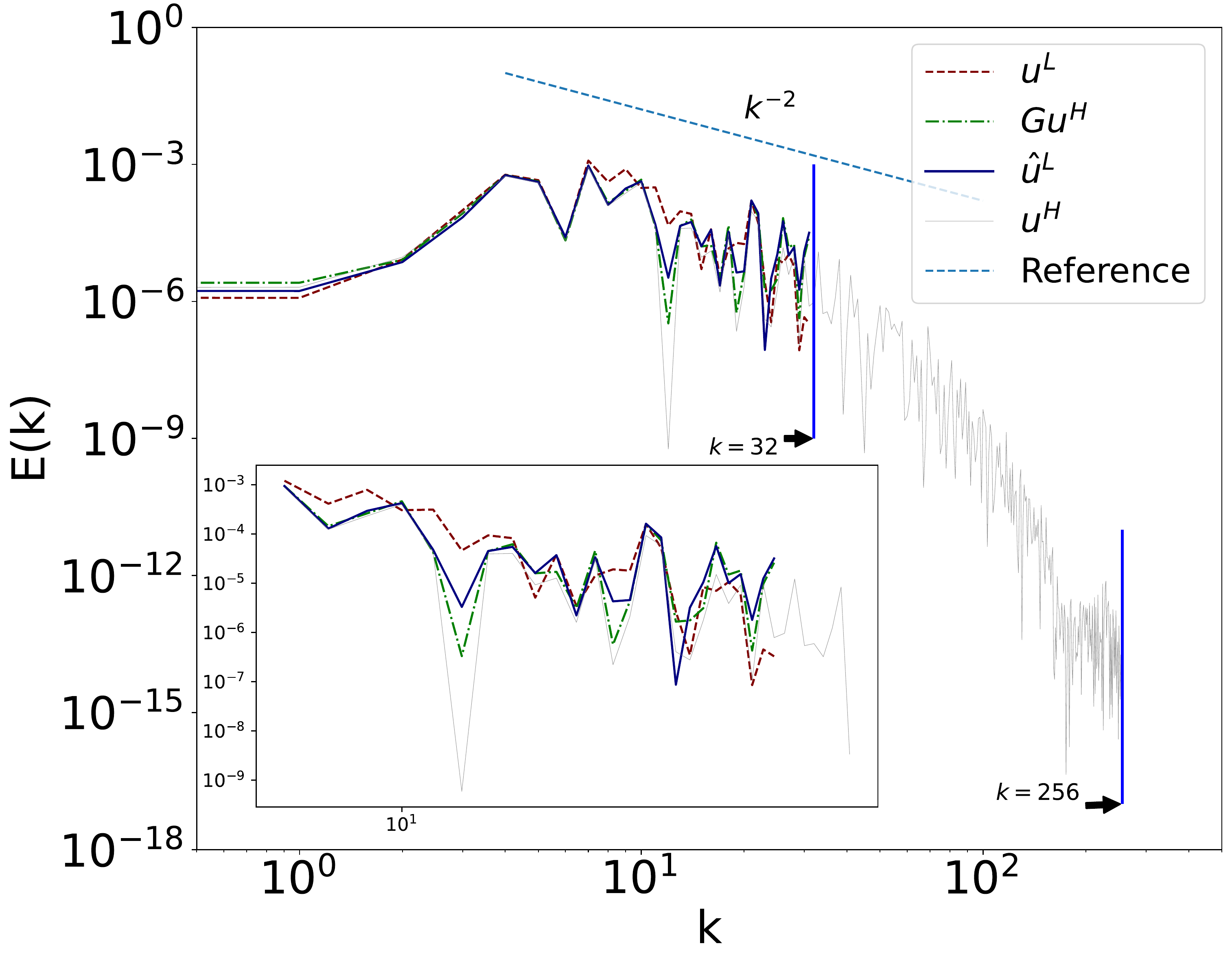}
}
  \caption{Kinetic energy spectra of $u^L$, $\hat{u}^L$, $Gu^H$, and $u^H$ with $\Nel=64$ at $t=0.5$ and (b) $t=1.0$. 
  }
  \figlab{burgers-spectra}
\end{figure}

  To see the consistency of 
 $\hat{u}^L$ prediction with respect to the timestep size, we ran two simulations with $\dt^L=10\times \dt^H$ and $\dt^L=50\times \dt^H$.
Figure \figref{Burgulence-predxt-diff} shows the difference of $u^L$ and $\hat{u}^L$ with respect to $Gu^H$ in a range between $-0.5$ and $0.5$. 
The absolute value of $u^L - Gu^H$ is $0.51$, and the 
 $\hat{u}^L - Gu^H$ counterpart is $0.06$ in both Figure \figref{Burgulence-predxt-diff}a and Figure \figref{Burgulence-predxt-diff}b. 
 The maximum $L_2$ errors of $u^L$ and $\hat{u}^L$ with respect to $Gu^H$ are $\norm{u^L - Gu^H}_{DG}=0.26$ and $\norm{\hat{u}^L - Gu^H}_{DG}=0.03$, respectively.
 We again observe that $\hat{u}^L$ is approximately nine times closer to $Gu^H$ than $u^L$. Also, the $\hat{u}^L$ prediction is consistent with respect to the timestep size. 

Next, we consider the polynomial equivalent of the neural network augmented approximation.
We integrate $u^L$, $\hat{u}^L$, $u_2^L$, and $u_3^L$ with $\dt=0.005$ for $t=[0,1]$ over
the mesh of $64$ elements. 
We measure the relative $L_2$ errors of $u^L$, $\hat{u}^L$, $u_2^L$, and $u^L_3$ by taking $Gu^H$ as the ground truth in Figure \figref{burgers-relerr-history}. For $t=[0,0.5]$, the augmented DG approximation $\hat{u}^L$ shows better accuracy than the others. 
After that, the second-order and the third-order DG approximations start to catch up with the accuracy level of $\hat{u}^L$. After $t=0.7$, the relative error of $u^L_3$ is smaller than that of $\hat{u}^L$.  After $t=0.8$, the relative error of $u^L_2$ is less than that of $\hat{u}^L$. 
When compared to the first-order DG approximation $u^L$, the enhanced DG approximation $\hat{u}^L$ still shows less relative error.


\begin{figure}[h!t!b!]
  \centering
  
\includegraphics[trim=0.0cm 0.5cm 0.2cm 0.2cm,clip=true,width=0.80\textwidth]{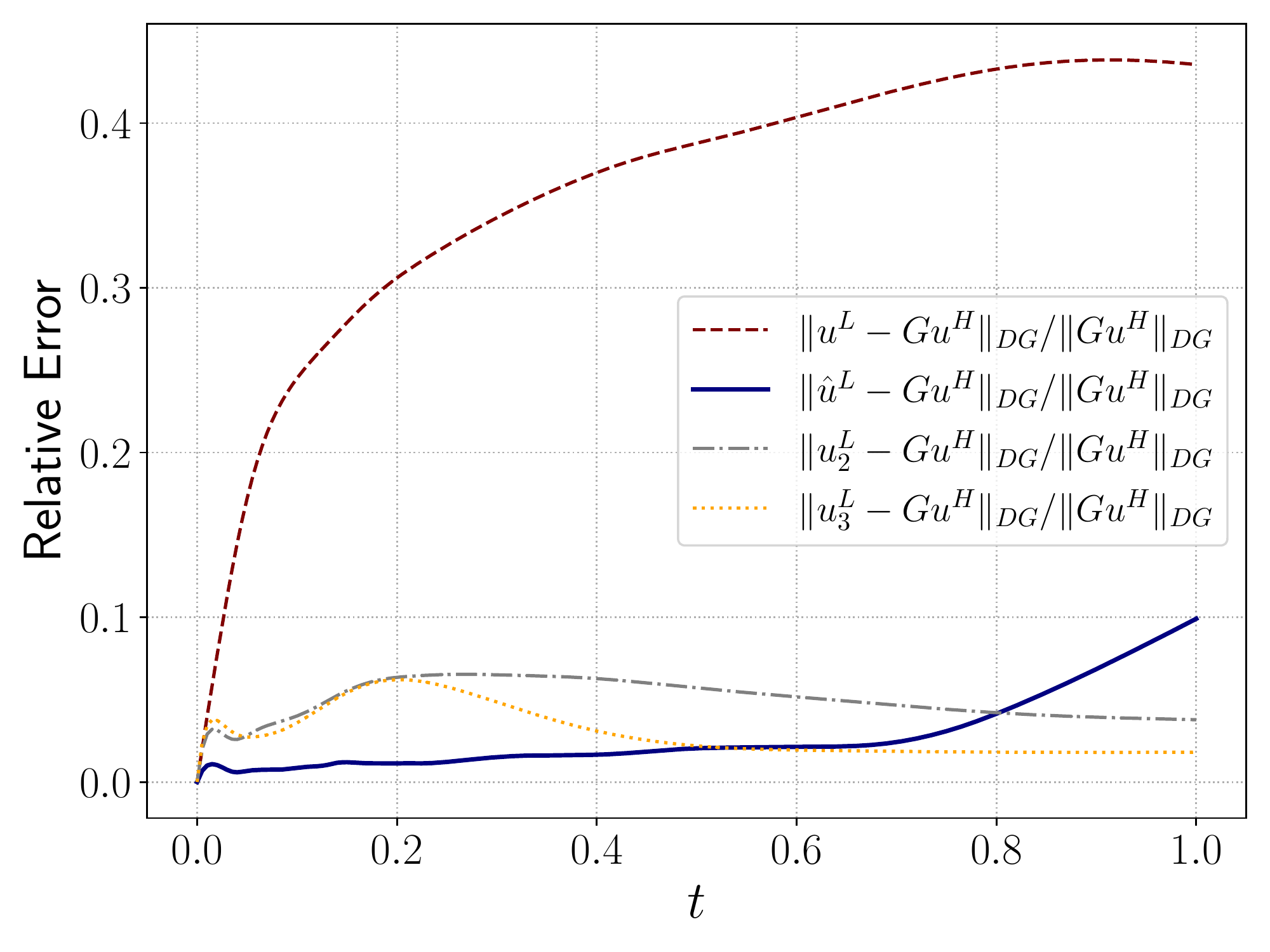}

  \caption{Burgulence: history of relative errors of $u^L$, $\hat{u}^L$, $u^L_2$, and $u^L_3$. The filtered solution $G u^H$ is considered as a reference solution.
  }
  \figlab{burgers-relerr-history}
\end{figure}

\begin{figure}[h!t!b!]
  \centering
   
\subfigure[$\dt^{L}=10\times \dt^H$]{  \includegraphics[trim=0.3cm 0.1cm 0.2cm 0.2cm,clip=true,width=0.4\textwidth]{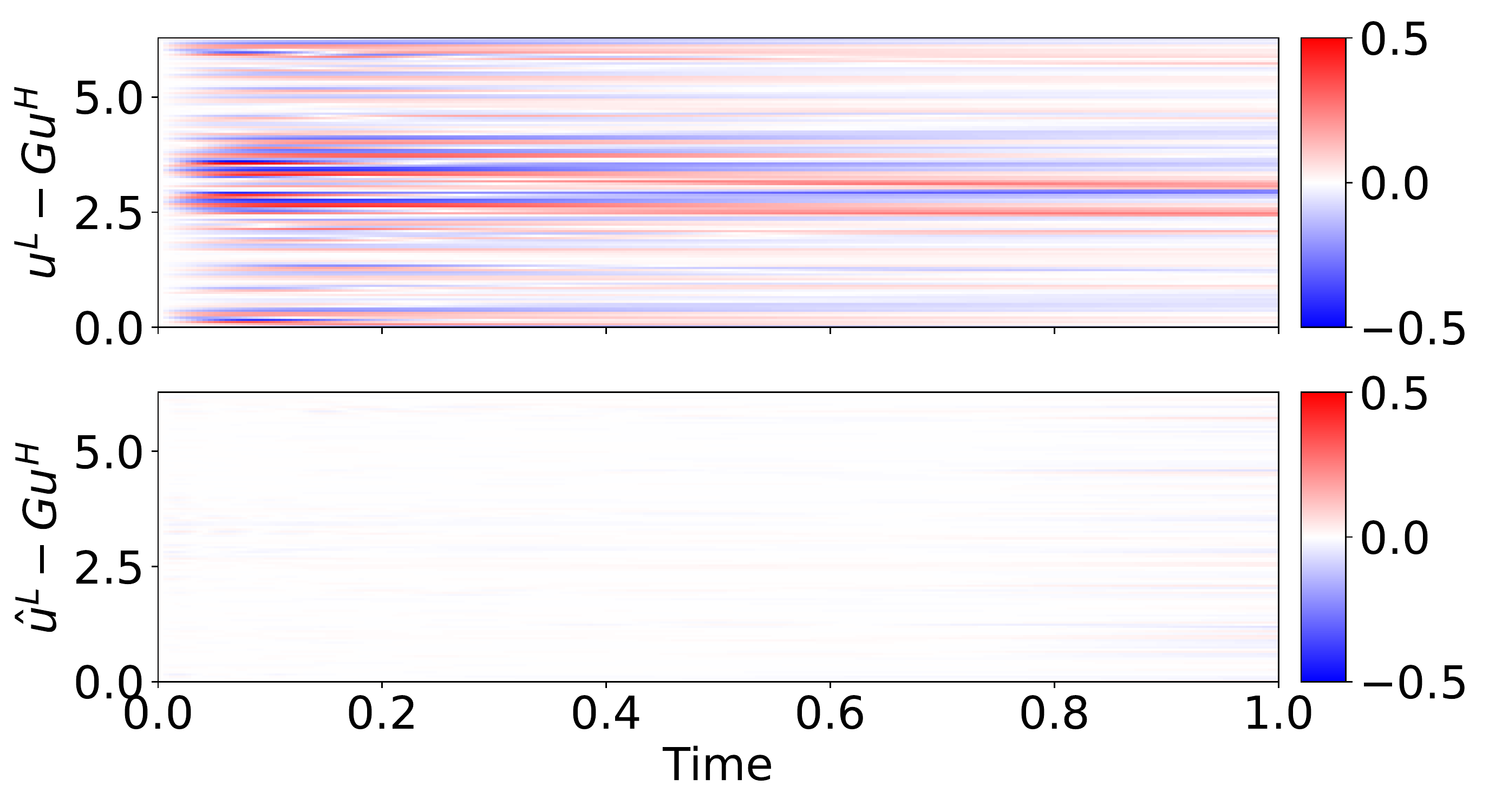}
}
\subfigure[$\dt^{L}=50\times \dt^H$]{  \includegraphics[trim=0.3cm 0.1cm 0.2cm 0.2cm,clip=true,width=0.4\textwidth]{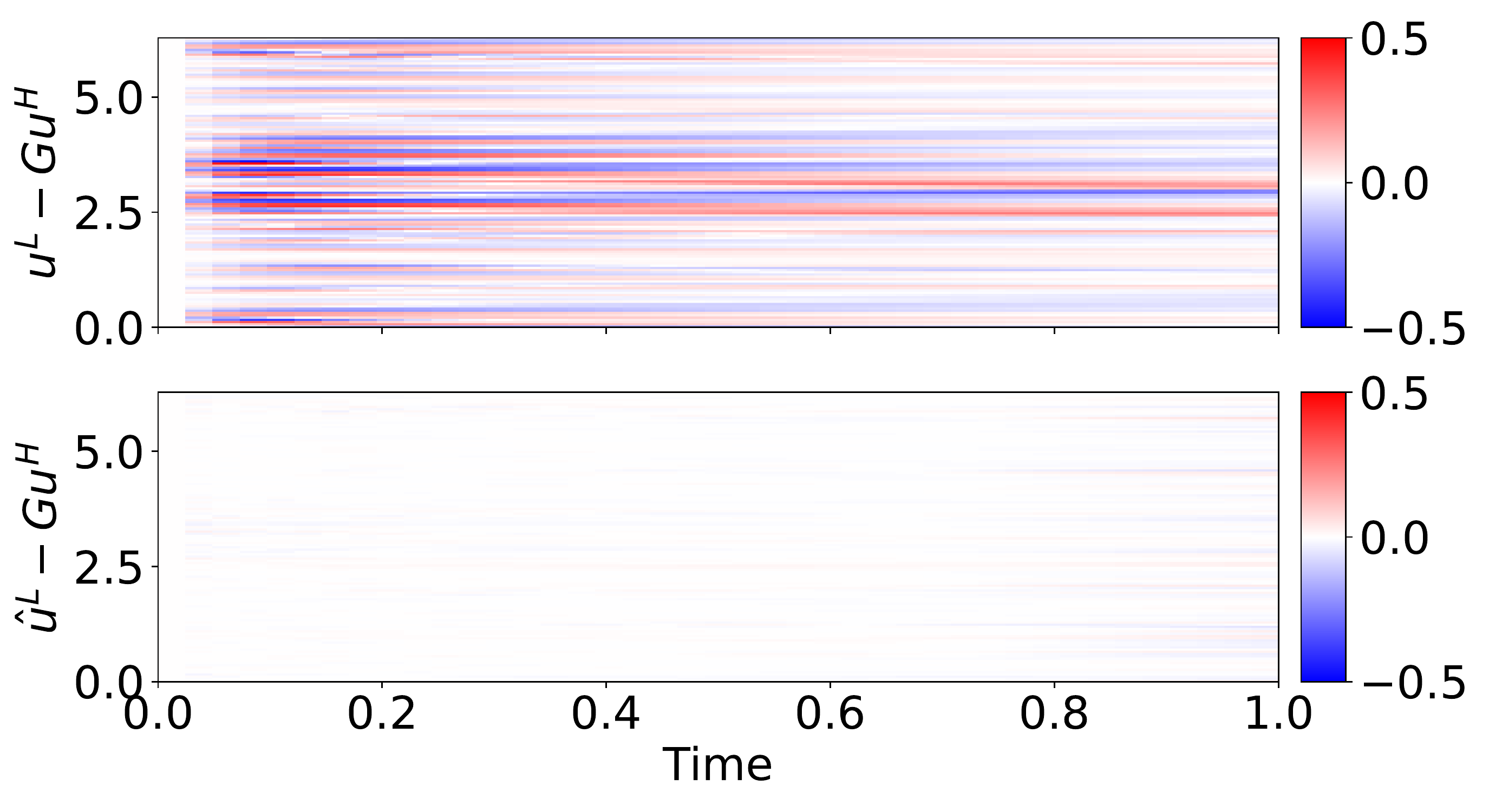}
}
  \caption{Burgulence: difference between $u^L$ and $G u^H$ (top)
  and the difference between $\hat{u}^L$ and $G u^{H}$ (bottom) with (a) $\dt^{L}=10\times \dt^H$ and (b) $\dt^{L}=50\times \dt^H$.
  }
  \figlab{Burgulence-predxt-diff}
\end{figure}


We also summarize the wall-clock time for Burgulence example in Table \tabref{vburgers-walltime}. 
We simulate this case ten times and report the average wall-clock time for the first and the second run.
We choose the stable timestep sizes of $\dt^H=0.001$ and $\dt^L=0.04$ so that increasing them by two causes a blowup. We denote the diffusion number by $Cr_d=\frac{\kappa \dt}{\dx^2}$. The Courant number and the diffusion number $(Cr_a,Cr_d)$ are ($0.23,0.21$) for $\dt^H=0.001$ and ($0.52,0.02$) for $\dt^L=0.04$. The timestep size of the high-order solution, $\dt^H$, is restricted by the diffusion number, whereas that of the low-order solution, $\dt^L$, is limited by the Courant number.  For $u^L_2$ and $u^L_3$, we take the stable timestep sizes of $\dt =0.018$ and $\dt=0.0125$, respectively. 
For the first run with JIT compiling time, all approximations are comparable, except that $\hat{u}^L$ is $25\%$ more expensive than the others.
In the second run, 
the wall clock of the low-order solution $u^L$ is $43$ times smaller than that of the high-order solution $u^H$. The neural network augmented low-order solution $\hat{u}^L$ is about $19$ times faster than the high-order solution $u^H$. The wall clock of $\hat{u}^L$ is 2.5 times larger than that of $u^L$, but comparable to that of $u^L_3$.

\begin{table}[t] 
    \caption{Viscous Burgers' equation: wall-clock time for the predictions of $u^H$, $u^L$, $\hat{u}^L$, $u^L_2$, and $u^L_3$.  }
    \tablab{vburgers-walltime} 
    \begin{center} 
    \begin{tabular}{*{1}{c}|*{1}{c}|*{1}{c}|*{1}{c}} 
    \hline 
      & $\dt $  & JIT Compile wc~[$ms$] & Simulation wc~[$ms$] \tabularnewline 
    \hline\hline 
    $u^H$       &  $0.001$ &     920 & 34.6  \tabularnewline
    $u^L$       &   $0.04$ &     919 &  0.8  \tabularnewline
    $\hat{u}^L$ &   $0.04$ &    1083 &  1.8  \tabularnewline
    $u^L_2$     & $0.018$  &     870 &  1.4  \tabularnewline
    $u^L_3$     & $0.0125$ &     858 &  1.9  \tabularnewline
    \hline\hline 
    \end{tabular} 
    \end{center} 
  \end{table} 

\section{Conclusions}
\seclab{Conclusion} 

In this work we present a methodology for learning subgrid-scale models based on neural ordinary differential equations. We utilize the NODE and partial knowledge to learn the underlying continuous 
source dynamics. The key idea is to augment the governing equations with a neural network source term and to train the neural network through the NODE. Our approach inherits the benefits of the NODE and can be applied to parameterize
subgrid scales, approximate coupling terms, and improve the efficiency of low-order solvers. 

We illustrate the proposed methodology on three model problems: the two-scale Lorenz 96 system, the linear convection-diffusion equation, and the viscous Burgers' equation. In the Lorenz 96 model, the slow and the fast variables are coupled to each other. In particular, the average of the fast variables contributes the dynamics of the slow variables. We replace the coupling component with the neural network source term. The prediction with the neural network coupling term shows good agreement with the true slow variable around $t=2\sim 3$ given a particular initial condition. This holds true even with a $10$ times larger timestep size. 
For the linear convection-diffusion equation, we include a neural network as a source term. Numerical examples demonstrate that adding a neural network source term considerably improves low-order DG approximation accuracy. For both the cases with $\kappa=10^{-4}$ and $\kappa=10^{-3}$, prediction with the neural network source term is ten times closer to the filtered high-order solutions than with the first-order DG approximation. For the case with $\kappa=10^{-4}$, the neural network augmented DG approximation is comparable to the projected third-order DG approximation. 
 For the case with $\kappa=10^{-3}$, the neural network augmented DG approximation enhances the first-order DG approximation but is not equivalent to the projected second-order DG approximation.
Our approach can be thought of as a generalization of the discrete corrective forcing approach, in which the first-order Euler method is employed to close the gap between the filtered and low-order solutions. 
Our approach learns continuous source dynamics and therefore can accommodate adaptive timestepping during training and prediction.
Indeed, the numerical examples confirm that our proposed approach is insensitive to changing the timestep size. In particular, the prediction with continuous corrective forcing term shows a convergent behavior to the filtered solution with decreasing timestep size. 
Moreover, we numerically demonstrate that our proposed approach performs better than the discrete corrective forcing approach in terms of accuracy, except for the case with $\dt=10^{-3}$ where the discrete corrective forcing was trained. The predictions with the neural network source term are ten times closer to the filtered high-order solutions than with the first-order DG approximation. Similarly, in the viscous Burgers' equation, the neural network source term enhances the low-order DG approximation accuracy. This improvement is also observed in the energy spectrum. The spectrum of the augmented low-order DG approximation is well overlapped with the spectrum of the filtered high-order solution. 
We also reported the wall-clock times for the low-order DG approximation, the augmented low-order DG approximation, and the high-order DG approximation. The low-order DG approximation with the neural network source term is about three to four times more expensive than the low-order DG approximation. Compared with the high-order DG approximation, augmenting the neural network source is $3$ and $19$ times more economical for the convection-diffusion model and the viscous Burgers' model, respectively. The wall-clock times for the enhanced DG solution are similar to those of the third-order DG solution in both the convection-diffusion and Burgulence examples.

Current work illustrates the potential applicability of employing a neural network source term with the aid of neural ordinary differential equations. However, our current approach has limitations. By introducing the neural network source term, the DG discretization loses the conservation property. Also, the neural network source term is global. This choice will definitely affect the parallel performance for large-scale applications. Furthermore, our current study is restricted to one-dimensional examples. Therefore, we will concentrate our future work on multidimensional extensions as well as other partial differential equations, such as Navier--Stokes equations. Exploring various network architectures, considering long-term stability and the conservation property, and devising local source operators are of interest.


\section*{Acknowledgments}
This material is based upon work supported by the U.S. Department of Energy, Office of Science, Office of Advanced Scientific Computing Research -- the Applied Mathematics Program -- and Office of Biological and Environmental Research, Scientific Discovery through Advanced Computing (SciDAC) program under Contract DE-AC02-06CH11357 through the Coupling Approaches for Next-Generation Architectures (CANGA) project and the FASTMath Institute. 
We also gratefully acknowledge the use of Theta, ThetaGPU, and Polaris resources in the Argonne Leadership Computing Facility, which is a DOE Office of Science User Facility supported under Contract DE-AC02-06CH11357.


\bibliographystyle{elsarticle-num}
\bibliography{main}

\end{document}